\documentclass{article}

 \usepackage[utf8]{inputenc}

\usepackage{amsmath,amsfonts,amssymb,amsthm}
\usepackage{bbm}

\usepackage[export]{adjustbox}
 \usepackage{setspace}
	
\usepackage{bm}
\usepackage{graphicx}
\usepackage{svg}
\usepackage{pgf,tikz,pgfplots}
\usepackage{subfiles}

\usepackage{hyperref}
\usepackage{aliascnt}

\hypersetup{%
  colorlinks=true,
  linkcolor=blue,
  citecolor=blue,
  linkbordercolor=red,
  pdfborderstyle={/S/U/W 1}
}

\usepackage{caption}
\usepackage{subcaption}
\usepackage[ruled,vlined]{algorithm2e}
\usepackage{epstopdf}

\tikzset{
point/.style={circle,draw=black,inner sep=0pt,minimum size=3pt}
}
\pgfplotsset{
    soldot/.style={color=blue,only marks,mark=*}
    }

    \usepackage{csvsimple}

\pgfplotsset{compat=1.14}
\usepackage{mathrsfs}
\usetikzlibrary{arrows,calc,patterns}

\renewcommand{\[}{\begin{equation}}
\renewcommand{\]}{\end{equation}}
\renewcommand{\O}{\Omega}

\newcommand{\T}{\mathcal{T}}

\renewcommand{\v}{\mathbf{v}}
\renewcommand{\u}{\mathbf{u}}

\newcommand{\ttau}{\boldsymbol{\tau}}

\newcommand{\hV}{\hat{\V}}
\newcommand{\hphi}{\hat{\phi}}

\newcommand{\hu}{\hat{u}}

\newcommand{\hGamma}{\hat{\Gamma}}

\newcommand{\hbeta}{\hat{\beta}}
\newcommand{\hK}{\hat{K}}
\renewcommand{\hbeta}{\hat{\beta}}

\newcommand{\n}{\mathbf{n}}

\renewcommand{\H}{\mathcal{H}}

\newcommand{\hPi}{\hat{\Pi}}

\newcommand{\tGamma}{\tilde{\Gamma}}

\usepackage{stmaryrd}
\newcommand{\bb}[1]{\left\llbracket#1\right\rrbracket}

\newcommand{\cc}[1]{\left\{\!\!\left\{#1\right\}\!\!\right\}}

\newcommand{\diam}{\operatorname{diam}}

\renewcommand{\c}{\mathbf{c}}
\newcommand{\x}{\mathbf{x}}
\renewcommand{\b}{\mathbf{b}}
\newcommand{\g}{\mathbf{g}}

\newcommand{\E}{\mathcal{E}}
\renewcommand{\L}{\mathscr{L}}

\newcommand{\norm}[1]{\left\lVert#1\right\rVert}

\newcommand{\N}{\mathcal{N}}

\renewcommand{\P}{\mathcal{P}}

\newcommand{\V}{\mathcal{V}}

\newcommand{\p}{\partial}

\newcommand{\ds}{\displaystyle}

\usepackage[utf8]{inputenc}

\usepackage{cases}

\newcommand{\btimes}{\wedge}

\newcommand{\comment}[1]{}
\newcommand{\draft}
{
\renewcommand{\comment}[1]{
\begin{center}
\fbox{\parbox{.8\textwidth}{
{\color{red}Comment:  }##1
}}
\end{center}
}
}

\usepackage{scalerel}

\usepackage{cancel}

\usepackage[top=2cm, bottom=2cm,left=2cm,right=2cm]{geometry}

\usepackage{enumitem}

\usepackage{cleveref}

\theoremstyle{definition}

\newtheorem{lemma}{Lemma}
\newtheorem{theorem}{Theorem}

\newtheorem{corollary}{Corollary}

\title{A High Order  Geometry Conforming Immersed Finite Element for  Elliptic Interface Problems }

\author{
  Slimane Adjerid\thanks{adjerids@vt.edu},
  Tao Lin\thanks{tlin@vt.edu}, and
   Haroun Meghaichi\thanks{haroun@vt.edu} \\ \\
  Department of Mathematics MC 0123,  Virginia Tech, 225 Stanger St. 460 McBryde Hall, \\ Blacksburg VA 24061 USA
}

\date{August 2023}

\newcommand{\commentout}[1]{{}} 
\draft

\begin{document}

\maketitle

\begin{abstract}
    We present a high order immersed finite element (IFE) method  for solving the elliptic interface problem with interface-independent meshes. The IFE functions developed here satisfy the interface conditions exactly and they have optimal approximation capabilities. The construction of this novel IFE space relies on a nonlinear transformation based on the Frenet-Serret frame of the interface to locally map it into a line segment, and this feature makes the process of constructing the IFE functions cost-effective and robust for any degree. This new class of immersed finite element functions is locally conforming with the usual weak form of the interface problem so that they can be employed in the standard  interior penalty discontinuous Galerkin scheme without additional penalties on the interface. Numerical examples are provided to showcase the convergence properties of the method under $h$ and $p$ refinements.
\end{abstract}


\section{Introduction}
In this paper, we introduce a new immersed finite element (IFE) approach for the second order elliptic interface problem \begin{subequations}\label{eqn:standard_problem_statement}

\begin{equation}
    \begin{cases}
        -\nabla\cdot (\beta \nabla u) = f,& \text{on } \Omega^+\cup\Omega^-,\\
        u_{\mid \partial \Omega} =g,
    \end{cases}
   \label{eqn:Model_problem_elliptic}
\end{equation}
where, without loss of generality, the domain $\Omega\subset \mathbb{R}^2$ is assumed to be a rectangle divided by an interface $\Gamma$ into two subdomains $\O^+$ and $\O^-$, and the diffusion coefficient $\beta$ is piecewise constant: $\beta|_{\O^{\pm}}= \beta^{\pm}>0$. In addition to \eqref{eqn:Model_problem_elliptic}, the solution $u$ is assumed to satisfy the following  jump conditions at the interface $\Gamma$
\begin{equation}
    \bb{u}_{\Gamma}=0, \label{eqn:continuity_cond}
\end{equation}
\begin{equation}
    \bb{\beta\frac{\p u}{\p\n} }_{\Gamma}=0, \label{eqn:normal_cond}
\end{equation}
\end{subequations}
where $\bb{v}_{\Gamma}=v^+|_{\Gamma} -v^-|_{\Gamma}$ with $v^{\pm}=v|_{\Omega^{\pm}}$ denotes the jump of $v$ across the interface $\Gamma$ and $\n$ is a normal vector of $\Gamma$. Also, we restrict our work to the case where $f$ is smooth around the interface, equivalently
\begin{equation}
    \bb{\beta\frac{\p^j\Delta u}{\p\n^j} }_{\Gamma}=0\qquad j=0,1,\dots, m-2. \label{eqn:laplacian_cond}
\end{equation}

%

The interface problem described by \eqref{eqn:standard_problem_statement} appears in numerous applications, such as topology optimization of heat conduction \cite{gersborg-hansen_topology_2006}, electrostatic levitation of dust particles \cite{wang_modeling_2008}, projection methods for Navier-Stokes equations \cite{chorin_numerical_1968}, water and oil reservoir simulation \cite{1983Ewingreservoir}, to name a few.

An essential idea of the IFE method is to use IFE functions on each interface element while standard finite element functions in terms of polynomials are employed on non-interface elements as usual. Instead of a polynomial, an IFE function in  the related literature is a piecewise polynomial, like the well-known Hsieh-Clough-Tocher macro $C^1$ element \cite{clough_finite_1966}, constructed according to the jump conditions for an interface problem. Generally speaking, an interface problem consists of two parts, the partial differential equation for the physics to be modeled (plus the usual boundary/initial conditions) in a domain
$\Omega$ and the jump conditions across an interface $\Gamma$ that partitions the modeling domain $\Omega$ into sub-domains. Obviously, the interface $\Gamma$ itself is a critical component across which the solutions to the PDE in the sub-domains are pieced together according to the relevant physics involved. Nevertheless, most, if not all, IFE methods in the literature use the information of the interface curve $\Gamma$ rather rudimentarily.
For example, only two intersection points of the interface $\Gamma$ and the edge of an interface element are used in the construction of IFE functions based on linear \cite{li_immersed_2004, li_new_2003} or bi-linear polynomials \cite{he_approximation_2008, lin_rectangular_2001} on that interface element. For constructing IFE functions based on higher degree polynomials, more points on the interface are used \cite{benromdhane2014} or the interface is employed in line integrals for weakly enforcing the jump conditions \cite{adjerid_higher_2018}. Even in more advanced techniques for constructing higher degree IFE functions such as the
least-squares method \cite{adjerid_high_2017} or Cauchy extension method \cite{guo_higher_2019}, the interface curve $\Gamma$ is still just used in integrals for approximately maintaining jump conditions. It is worth noting that, along a curve, its tangent vector field, normal vector field, and curvature are fundamental differential geometry characteristics \cite{gray_2006, ONeill_DiffGeom_2010}. However, as commented above, in the existing literature on the construction of IFE functions, these properties have not been explicitly utilized yet. {In addition to the IFE methods mentioned above, other numerical methods have been developped and analyzed for the interface problem \eqref{eqn:standard_problem_statement} such as the CutFEM method \cite{burman_cutfem_15} and the immersed interface method \cite{leveque_immersed_1994} among others.  } \\

In this article, we propose a new method based on the differential geometry of the interface $\Gamma$ for constructing IFE functions to solve the elliptic interface problem described by \eqref{eqn:standard_problem_statement}. Under the assumption that the interface $\Gamma$ has a regular parametrization so that it has nonzero tangential vector \cite{Pressley_DiffGeom_2010}, the tangential vector and normal vector at a point on the curve $\Gamma$ form the Frenet frame (coordinate system) \cite{gray_2006, ONeill_DiffGeom_2010} which can be used to represent the points in the neighborhood of this interface point. This means that a function in terms of the usual Cartesian coordinates $x$-$y$ can be expressed as a function in terms of the Frenet coordinates $\eta$-$\xi$ which are intrinsically linked to the curve $\Gamma$ in this neighbourhood. This Frenet transformation naturally admits the fundamental differential geometry properties such as the tangent vector, the normal vector, and the curvature into the formulas of jump conditions which can be subsequently employed in the construction of IFE functions across the interface.

To be  {more specific},
we introduce a fictitious element $K_F$ for each interface element $K$ in a mesh $\mathcal{T}_h$ of the solution domain $\Omega$. This fictitious element is a local Frenet tube of the interface $\Gamma$ that tightly covers the interface element, and it is a curved trapezoid formed by two curves parallel to the interface curve $\Gamma$ and two straight sides normal to the interface $\Gamma$, see the illustration in \autoref{fig:intf_fic_elements}. Through the Frenet transformation, this fictitious element is the image of a rectangle $\hat{K}_F$ in the Frenet coordinate system and the interface section $\Gamma_{K_F} = \Gamma \cap K_F$ is the image of the line segment in the $\xi$-axis inside the rectangle $\hat{K}_F$. These are key geometric ideas for the method proposed here to construct IFE functions on interface elements, which brings in numerous advantages over methods in the existing literature. This allows us to develop IFE functions with polynomials in the fictitious element $\hat{K}_F$ in the Frenet coordinate system across a simple interface that is a straight line, and then obtain IFE functions on the interface element through the Frenet transformation. Of course, similar to the well-known isoparametric IFE functions \cite{D.Braess, ciarlet_78}, IFE functions constructed in the interface element $K$ in this way are not polynomials anymore because of the nonlinearity of the Frenet transformation. Nevertheless, this key idea greatly facilitates the construction of IFE functions with higher degree polynomials, and the proposed method have the following distinct features:

\begin{itemize}[label=\textbullet] 
\item
With the simple interface geometry in the Frenet fictitious element $\hat{K}_F$, the existence of a basis of IFE functions based on general $Q^m$ polynomials is established regardless of the complexity of the interface $\Gamma$. A basis of IFE functions can be constructed with much less computations and the involved linear system in the construction is well conditioned regardless of the interface location and the magnitude of the contrast in the coefficient $\beta$.

\item
Construction methods in the literature produce higher degree IFE functions that can satisfy jump conditions specified in \eqref{eqn:continuity_cond} and \eqref{eqn:normal_cond} {\em only approximately}. In contrast, the IFE functions produced by the proposed method will satisfy the jump conditions \eqref{eqn:continuity_cond} and \eqref{eqn:normal_cond} {\em precisely}, and for higher degree polynomials, the majority of the basis functions
in the local IFE space on each interface element can even satisfy the extended jump conditions \eqref{eqn:laplacian_cond} {\em precisely}.

\item
Being able to satisfy interface jump conditions \eqref{eqn:continuity_cond} and \eqref{eqn:normal_cond} makes the local IFE space on each interface element conforming to the Sobolev space employed in the weak form of the interface problem described by \eqref{eqn:standard_problem_statement}. This has great potential for alleviating challenges in the error estimation of numerical methods based on the IFE space proposed here.

\item
The proposed IFE space has the optimal approximation capability with respect to the underlying $Q^m$ polynomials and the error of the approximation is independent of the coefficients $\beta^\pm$ and relative position of the interface to the mesh.

\end{itemize}
We observe that none of the methods in the literature for constructing IFE functions have the first three features listed above.

The rest of this article is organized as follows, in \autoref{sec:Prelim}, we introduce the notation and spaces used in the rest of the paper. In \autoref{sec:Frenet_coordinates}, we describe the local coordinate system and the Frenet transformation. In \autoref{sec:Frenet_space}, we discuss the construction of the local IFE space and its properties. In \autoref{sec:Approximation}, we prove that the proposed IFE space has optimal approximation capabilities. In \autoref{sec:the_ipdg_method}, we present an immersed DG method to solve the elliptic interface problem, as well as some details for the implementation of the Frenet transformation and the construction of the local matrices. At the end, in \autoref{sec:numerical_examples}, we provide some numerical examples to illustrate the features of our proposed method.

\section{Preliminaries}\label{sec:Prelim}
We assume throughout this manuscript, without loss of generality, that $\Omega\subset \mathbb{R}^2$ is a rectangular domain  which is split by an interface curve $\Gamma$ into two sub-domains
$\Omega^s, s = \pm$ such that $\p\Omega^+\cap\p\Omega^{-}=\Gamma$. For a given measurable set $\tilde{\O}$, we denote the standard Sobolev space on this set by $H^{s}(\tilde{\O})$ with the norm $\norm{\cdot}_{s,\tilde{\O}}$.  We also use $(\cdot,\cdot)_{\tilde{\O}}$ and $\langle \cdot ,\cdot \rangle_{\p \tilde{\O}} $ to denote the standard $L^2$ inner product on $\tilde{\O}$ and $\p \tilde{\O}$, respectively.
If $\tilde{\O}\subset \O$ intersects $\O^-$ and $\O^+$, we let $\tilde{\O}^{s}=\tilde{\O}\cap \O^s, s = -, +$ and we define the immersed Sobolev space

\begin{equation}
\H^s(\tilde{\O},\Gamma;\beta)
=\left\{u:\tilde{\O}\to \mathbb{R}\mid u_{\tilde{\O}^\pm}\in H^{s}(\tilde{\O}^\pm),\ \bb{u}_{\Gamma\cap \tilde{\O}}=\bb{\beta \nabla u\cdot \n}_{\Gamma\cap \tilde{\O}} =0\right\},
\end{equation}
where we assume $s>\frac{3}{2}$ for the involved quantities to be well-defined. We equip $\H^s(\tilde{\O},\Gamma;\beta)$ with the broken Sobolev norm $\norm{\cdot}_{s,\tilde{\O}} =\sqrt{\norm{\cdot}_{s,\tilde{\O}^+}^2+\norm{\cdot}_{s,\tilde{\O}^-}^2 }$. The related semi-norms for these spaces are defined similarly.

For the solution domain $\Omega$, we consider an interface-independent uniform mesh $\T_h$ of rectangles/squares with {diameter} $h$. However, we note that all of our discussions in this paper apply to non-uniform meshes with quadrilateral or triangular elements.  We call an element $K\in \T_h$ an interface element if $\mathring{K}\cap \Gamma \neq \varnothing$, where $\mathring{K}$ is the interior of $K$; otherwise, we call $K$ a non-interface element. We use $\T_h^i$ and $\T_h^n$ to denote the set of interface elements and the set of non-interface elements, respectively. In addition, we denote the set of edges, interior edges and boundary edges of $\T_h$ by $\E_h,~\E_h^\circ$ and $\E_h^b$. In the mesh $\T_h$, we define the broken immersed Sobolev space
\begin{equation}
\H^s(\T_h,\Gamma;\beta) = \left\{
u:\Omega\to \mathbb{R}\mid u|_{K} \in H^s(K) \text{ if } K\in\T_h^n;\text{otherwise } u|_{K}\in \H^s(K,\Gamma;\beta)\right\}.
\label{eqn:def_broken_immersed_sob}
\end{equation}

From now on, we assume that the solution $u$ to the problem  \eqref{eqn:standard_problem_statement} is in $\H^{s}(\O,\Gamma;\beta)\subset \H^{s}(\T_h,\Gamma;\beta) $, where $s>\frac{3}{2}$. This allows us to characterize $u\in \H^{s}(\O,\Gamma;\beta)$ as the solution to
\begin{subequations}
\begin{equation} a_h(u,v)=L_h(v),\qquad \forall v\in \H^{s}(\T_h,\Gamma;\beta),\label{eqn:compact_weak_form} \end{equation}
where $a_h:\H^{s}(\T_h,\Gamma;\beta)\times \H^{s}(\T_h,\Gamma;\beta)\to \mathbb{R} $ is the following symmetric bilinear form
\begin{eqnarray} \label{eqn:def_of_a_h}
&& a_h(u,v) =\sum_{K\in \T_h} \left(\beta\nabla u,\nabla v\right)_{K} \nonumber \\
&& \hspace{0.6in} -\sum_{e\in \E_h} \left(\left\langle \bb{\beta\nabla u\cdot \n_e}_e,\cc{v}_e\right\rangle_{e} +
\left\langle\cc{u}_e, \bb{\beta\nabla v\cdot \n_e}_e\right\rangle_{e} -
\frac{{\sigma_e}}{h}\left\langle\bb{u}_e, \bb{v}_e\right\rangle_{e}
\right),\label{eqn:def_a_h}
\end{eqnarray}
and $L_h: \H^{s}(\T_h,\Gamma;\beta)\to \mathbb{R}$ is the linear form defined as follows

\begin{equation}
    L_h(v)  = (f,v)_{\O} + \sum_{e\in\E_h^b} \left\langle -\beta \nabla v\cdot \n_e +\frac{{\sigma_e}}{h}v,g\right\rangle_e.
    \label{eqn:def_of_L_h}
\end{equation}
\end{subequations}
Here, $\n_e,\cc{\cdot}_e,\bb{\cdot}_e$ denote the normal vector to the edge $e$, the average of a function across an edge $e$, and the jump across an edge $e$, which are standard quantities used for DG methods \cite{riviere_discontinuous_2008}.
It is worth mentioning here that in the case where $\beta$ is constant, the formulation \eqref{eqn:compact_weak_form} is related to the classical symmetric interior penalty Galerkin (SIPG) method  \cite{riviere_discontinuous_2008}.

\section{The Frenet coordinates for interface elements} \label{sec:Frenet_coordinates}

In our discussions from now on, without loss of generality,
we  assume that the interface $\Gamma$ is a {simple} connected $C^2$ curve described by the following function:
\begin{eqnarray*}
\g(\xi) = (g_1(\xi), g_2(\xi))^T: [\xi_s, \xi_e] \rightarrow \Omega,
\end{eqnarray*}
but a $C^3$ smoothness of $\g(\xi)$ will be introduced when the situation mandates. We also assume that the parametrization $\g$ of $\Gamma$ is regular in the sense that $\g'(\xi) \not = {\bf 0}$ for all $\xi \in [\xi_s, \xi_e]$
\cite{Pressley_DiffGeom_2010,Milliam_Parker_1977}. This assumption allows us to employ fundamental differential geometry quantities of a curve in discussions/computations involving the interface $\Gamma$.

Specifically, we will use the following components in the Frenet apparatus \cite{Milliam_Parker_1977} of
the curve $\Gamma$: the unit tangent vector $\ttau(\xi)$ at a point $\g(\xi) \in \Gamma$ given by
\begin{eqnarray}
\ttau(\xi) = \frac{1}{\norm{\g'(\xi)}} \g'(\xi), \label{eq:tangent_vec}
\end{eqnarray}
the unit normal vector at $\g(\xi) \in \Gamma$ given by
\begin{equation}
\n(\xi)= Q\ttau(\xi),\qquad Q=\begin{bmatrix}0&1\\ -1& 0\end{bmatrix}, \label{eqn:def_of_Q}
\end{equation}
and the signed curvature $\kappa(\xi)$ at $\g(\xi) \in \Gamma$ defined by
\begin{eqnarray}
\kappa(\xi) = \frac{1}{\norm{\g'(\xi)}}\big(\g'(\xi) \wedge \g''(\xi)\big), ~~\text{where}~~\u\btimes \v= \u^T Q\v.
\label{eqn:curvature}
\end{eqnarray}

\begin{figure}[ht]
\centerline{
\includegraphics[scale=0.75]{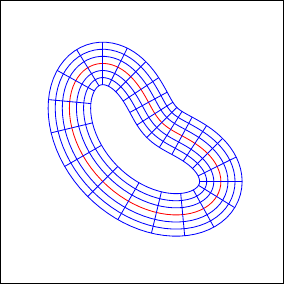}
}
        \caption{Curves parallel to $\Gamma$ and lines normal to $\Gamma$ form a curvilinear coordinate system in the vicinity of the interface $\Gamma$. The innermost curve, parallel to $\Gamma$, and the outermost curve, also parallel to $\Gamma$, enclose a region known as a tubular neighborhood of $\Gamma$.}
\label{fig:FerentRibbon_Curvilinear}
\end{figure}

In the vicinity of the interface $\Gamma$, there are curves parallel to $\Gamma$ in the sense that there exists a one-to-one mapping between each pair of these curves and $\Gamma$ such that corresponding points are equally distant and the tangents at corresponding points are parallel \cite{Lane_1940}. Furthermore, in terms of the Frenet apparatus, a curve parallel to $\Gamma$ with an offset distance $\eta$ from $\Gamma$ has the following parametric form \cite{Toponogov_2006}:
\begin{eqnarray}
\x(\eta, \xi) = \begin{bmatrix}
x(\eta, \xi) \\
y(\eta, \xi)
\end{bmatrix} = P_\Gamma(\eta, \xi) = \g(\xi) + \eta \n(\xi), ~\xi \in [\xi_s, \xi_e]. \label{eq:paralle_curve}
\end{eqnarray}
Straight lines normal to $\Gamma$ and curves parallel to $\Gamma$ form a curvilinear coordinate system in a neighborhood of $\Gamma$, see an illustration in \autoref{fig:FerentRibbon_Curvilinear}. Locally in this curvilinear system, the line normal to $\Gamma$ acts as the horizontal axis, and the curve perpendicular to this normal vector is the vertical axis. Since the direction along the curve is the tangent vector of the curve, this curvilinear coordinate system is determined by the two fundamental components of the Frenet apparatus: the normal vector and the tangent vector of the curve $\Gamma$. This suggests to use the Frenet coordinates $(\eta, \xi)$ to represent a point $(x, y)$ in Cartesian coordinates in this neighborhood, and we call the function given in \eqref{eq:paralle_curve} the Frenet transformation from the Frenet coordinates $(\eta, \xi)$ to the Cartesian coordinates. Since the Frenet coordinates intrinsically have fundamental features of the interface curve $\Gamma$, using the Frenet coordinates is advantageous in an IFE method that relies on IFE functions constructed according to jump conditions across the interface $\Gamma$. We note that
transformations based on the interface have been used in publications \cite{benromdhane2014, benromdhane2015, guzman2016, lehrenfeld_high_2016} for dealing with interface problems. However, to the best of the authors' knowledge,
this article is the first to use \eqref{eq:paralle_curve} extensively in constructing IFE functions with higher degree polynomials for solving interface problems.

In principle, we only need to consider the usage of Frenet coordinates locally in interface elements because an IFE method always uses standard polynomial type finite element functions on all the non-interface elements which have nothing to do with the interface. Without loss of generality, we assume that the interface curve
$\Gamma$ has a tubular neighborhood with a half-band width $\epsilon > 0$
denoted by
\begin{eqnarray}
N_\Gamma(\epsilon) = P_\Gamma([-\epsilon, \epsilon] \times [\xi_s, \xi_e]) \label{eq:tubular}
\end{eqnarray}
that has the following properties \cite{Abate_Tovena_2012,federer1959curvature}: {\it (i)} for every point $(x(\eta, \xi), y(\eta, \xi)) \in N_\Gamma(\epsilon)$, the Jacobian $\hat{D}P_\Gamma(\eta,\xi)$ of the Frenet transformation at $(\eta, \xi)$ is nonsingular; {\it (ii)} the lines perpendicular to $\Gamma$ passing any two distinct points on $\Gamma$ will not intersect within this tubular neighborhood $N_\Gamma(\epsilon)$; and {\it (iii)}, importantly,
the mapping $P_\Gamma: [-\epsilon, \epsilon]\times [\xi_s, \xi_e] \rightarrow N_\Gamma(\epsilon)$ is one-to-one and onto.
Because $P_\Gamma: [-\epsilon, \epsilon]\times [\xi_s, \xi_e] \rightarrow N_\Gamma(\epsilon)$ is a one-to-one and onto mapping, it has the inverse $R_\Gamma: N_\Gamma(\epsilon) \rightarrow [-\epsilon, \epsilon]\times [\xi_s, \xi_e]$ such that
\begin{eqnarray}
\begin{bmatrix}
\eta \\
\xi
\end{bmatrix} = \begin{bmatrix}
\eta(x, y) \\
\xi(x, y)
\end{bmatrix} = R_\Gamma(x, y)  = P_\Gamma^{-1}(x, y) \in [-\epsilon, \epsilon]\times [\xi_s, \xi_e],~~\forall \x \in N_\Gamma(\epsilon). \label{eq:P_inv}
\end{eqnarray}

\begin{figure}[ht]
\centerline{
\includegraphics[scale=0.75]{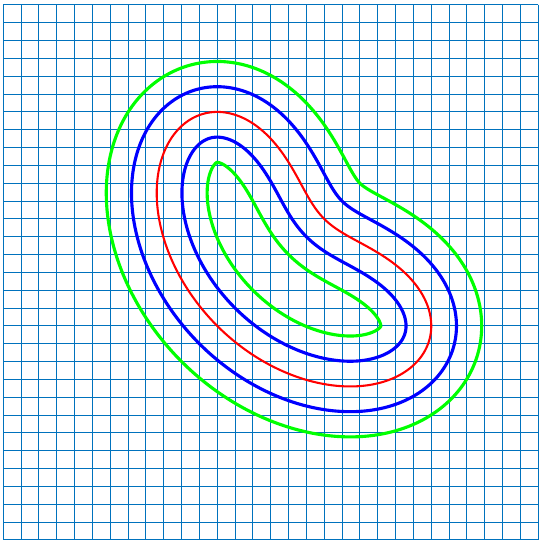}
}
        \caption{The two green parallel curves form an $\epsilon$ tubular neighborhood $N_\Gamma(\epsilon)$ of
        the interface curve $\Gamma$ in red color, the two blue parallel curves form the $h$ tubular neighborhood $N_\Gamma(h)$.}
\label{fig:FrenetRibbon}
\end{figure}

\begin{figure}[ht]
\centerline{
\includegraphics[scale=1]{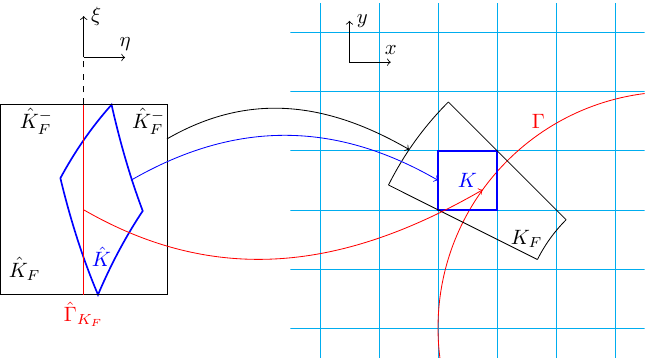}
\hspace{0.4in}
}
        \caption{The plot on the right is for a square interface element $K$ and the associated fictitious element $K_F$ which is the black trapezoid with two bases parallel to the interface curve $\Gamma$ in red color. The black rectangle in the plot on the left is the related fictitious element $\hK_F$ in the $\eta$-$\xi$ plane such that $K_F = P_\Gamma(\hK_F)$.}
\label{fig:intf_fic_elements}
\end{figure}

Obviously, for a mesh $\mathcal{T}_h$ with $h$ as the maximum of diameters of all elements, all the interface elements are inside the tubular neighborhood $N_\Gamma(h)$ of $\Gamma$ because every interface element intersects  $\Gamma$. Hence, we can assume that the mesh size $h$ is small enough such that $N_\Gamma(h) \subset N_\Gamma(\epsilon)$ {or simply $h<\varepsilon$}, see an illustration in \autoref{fig:FrenetRibbon}. For each interface element $K = \square A_1A_2A_3A_4 \in \mathcal{T}_h^i$, we can use a section of this tubular neighborhood to form a fictitious element $K_F$ that covers $K$. One way to construct this fictitious element is to note that for each vertex $A_i, 1 \leq i \leq 4$ of $K$, with $A_i \in N_\Gamma(\epsilon)$, there exists a unique $\xi_i \in [\xi_s, \xi_e]$ such that \cite{Abate_Tovena_2012}
\begin{eqnarray*}
\norm{A_i - \g(\xi_i)} = \text{dist}(\Gamma, A_i).
\end{eqnarray*}
We use $\xi_i, 1 \leq i \leq 4$ to determine two parameters $a_K, b_K \in [\xi_s, \xi_e]$ such that
\begin{eqnarray*}
{ a_K= \min(\xi_1, \xi_2, \xi_3, \xi_4),\qquad b_K= \max(\xi_1, \xi_2, \xi_3, \xi_4).}
\end{eqnarray*}
Then we form a rectangle $\hK_F = [-h, h] \times [a_K, b_K]$ and define the fictitious element $K_F$ of the interface element $K$ as
\begin{eqnarray}
K_F = P_\Gamma(\hK_F) = P_\Gamma([-h, h] \times [a_K, b_K]). \label{eq:fic_elem}
\end{eqnarray}
In the following discussions, for every interface element $K \in \mathcal{T}_h^i$, we will refer to $\hK = R_\Gamma(K)$ as the associated Frenet interface element of $K$ and
refer to $\hK_F$ as the Frenet fictitious element of $K$, see illustrations in \autoref{fig:intf_fic_elements}.

Geometrically, the fictitious element $K_F$ is formed by two straight lines aligned with the normal
vectors of $\Gamma$ at $a_K$ and $b_K$, respectively, and the two curves parallel to $\Gamma$ with the offset distance $-h$ and $h$, see the illustration in \autoref{fig:intf_fic_elements}, \textit{i.e.}, $K_F$ is a trapezoid with two curved bases parallel to the interface $\Gamma$ and two legs normal to the interface $\Gamma$. A critical feature of this fictitious interface element $K_F$ is that the interface curve segment $\Gamma_{K_F} = \Gamma \cap K_F$ is transformed into the Frenet fictitious interface $\hGamma_{K_F}$ which is a line segment
along the $\xi$ axis inside the corresponding Frenet fictitious element $\hK_F$. We will show later that this feature greatly simplifies the construction of IFE functions in terms of polynomials of $\eta$ and $\xi$ on $\hK_F$, and IFE functions on the interface element $K$ are subsequently formed by these IFE polynomials in $\eta$ and $\xi$ through the Frenet transformation \eqref{eq:paralle_curve}.

We now express a few basic calculus operators in terms of the Frenet coordinates $\eta$-$\xi$ which are essential in our discussions. To state the formulas clearly, we will use $\nabla, \nabla\cdot, D$ and $\Delta$ to denote the gradient, divergence, Jacobian matrix and Laplacian with respect to the Cartesian coordinates $\x=(x, y)$, respectively.
For the sake of conciseness and if there are no specific clarifications, in all the formulas from now on, we assume that the Frenet coordinates $(\eta, \xi)$ of a point inside the $\epsilon$-tube $N_\Gamma(\epsilon)$ of the interface curve $\Gamma$ are always associated with the Cartesian coordinate $\x = (x, y)$ of this point, {\it i.e.}, whenever
both $(\eta, \xi)$ and $\x = (x, y)$ appear in a formula, we will concisely use $(\eta, \xi)$ instead of $(\eta(\x), \xi(\x))$ whose dependence on $\x$ is determined by
\eqref{eq:P_inv}.

We start from recalling the well-known Frenet-Serret formulas \cite{gray_2006} for the tangent vector and the normal vector:
\begin{equation}
    \ttau'(\xi)=-\kappa(\xi)\norm{\g'(\xi)}\n(\xi),\qquad \n'(\xi)=\kappa(\xi)\norm{\g'(\xi)}\ttau(\xi).
    \label{eqn:frenet-serret}
\end{equation}
By \eqref{eqn:frenet-serret} we can show that the Jacobian of the Frenet transformation $P_\Gamma(\eta, \xi)$ is
\begin{equation}
    \hat{D}P_\Gamma(\eta,\xi)= \begin{bmatrix} \n(\xi), ~~\g'(\xi)+\eta\n'(\xi)\end{bmatrix}
    =\begin{bmatrix} \n(\xi), ~~\norm{\g'(\xi)}\left(1+\eta \kappa(\xi) \right)\ttau(\xi)\end{bmatrix},
    \label{eqn:jacobian_of_PK}
\end{equation}
where $\hat{D}P_{\Gamma}$ is the Jacobian matrix of $P_{\Gamma}$ with respect to $(\eta,\xi)$. By \eqref{eqn:jacobian_of_PK} and the inverse function theorem, we have the following formula for the Jacobian of  the inverse Frenet transformation $R_\Gamma(x, y)$:
\begin{equation}
    DR_\Gamma(\x)=\begin{bmatrix}
        \n(\xi)^T\\ \norm{\g'(\xi)}^{-1}\psi(\eta,\xi)\ttau(\xi)^T
    \end{bmatrix}, ~~\psi(\eta,\xi)= (1+\eta \kappa(\xi))^{-1}, ~~\forall \x \in N_\Gamma(\epsilon).
    \label{eqn:jacobian_of_RK}
\end{equation}
By the definition of $DR_\Gamma(\x)$ and \eqref{eqn:jacobian_of_RK}, we have
\begin{equation}
    \nabla \eta(\x) =\n(\xi),\qquad \nabla \xi(\x) = \norm{\g'(\xi)}^{-1}\psi(\eta,\xi)\ttau(\xi), ~~\forall \x \in N_\Gamma(\epsilon). \label{eqn:grad_eta_xi}
\end{equation}
Every function $u: N_\Gamma(\epsilon) \rightarrow \mathbb{R}$ is associated  with the following function of the Frenet coordinates:
\begin{eqnarray}
\hu = u\circ P_\Gamma: [-\epsilon, \epsilon]\times [\xi_s, \xi_e] \to \mathbb{R}. \label{eq:uhat}
\end{eqnarray}
Then, by the chain rule, we have the following formula for the gradient of $u(\x)$:
\begin{equation}
    \nabla u(\x) = \hu_{\eta}(\eta,\xi) \n(\xi) + \norm{\g'(\xi)}^{-1}\psi(\eta,\xi)\hu_{\xi}(\eta,\xi)\ttau(\xi), ~~\forall \x \in N_\Gamma(\epsilon).
    \label{eqn:grad_u}
\end{equation}
Consequently, we have
\begin{equation}
    \frac{\p u(\x)}{\p \n}=  \hu_{\eta}(\eta, {0}), ~~\forall \x \in \Gamma.
    \label{eqn:normal_to_eta}
\end{equation}
For the divergence of vector functions, we first can use \eqref{eqn:frenet-serret} and
\eqref{eqn:grad_eta_xi} to show that
\begin{equation}
    \nabla\cdot \n(\xi(\x))= \nabla\xi(\x) \cdot \n'(\xi)=\kappa(\xi)\psi(\eta,\xi), ~~\forall \x \in N_\Gamma(\epsilon),
    \label{eqn:div_eta}
\end{equation}
and
\begin{equation}
    \nabla\cdot \ttau(\xi(\x))= \nabla\xi(\x) \cdot \ttau'(\xi)=0, ~~\forall \x \in N_\Gamma(\epsilon).
    \label{eqn:div_xi}
\end{equation}
Given a generic vector function $\mathbf{f}: N_\Gamma(\epsilon) \to \mathbb{R}^2$, we can express it as
\begin{eqnarray*}
\mathbf{f}(\x) = f^{\n}(\x)\n(\xi) +f^{\ttau}(\x)\ttau(\xi),~~\forall \x \in N_\Gamma(\epsilon).
\end{eqnarray*}
Then, by \eqref{eqn:grad_u}, \eqref{eqn:div_eta}, and \eqref{eqn:div_xi}, we can show that
\begin{eqnarray}
   \nabla \cdot \mathbf{f}(\x) &=& \hat{f}^{\n}_{\eta}(\eta,\xi) + \kappa(\xi)\psi(\eta,\xi) \hat{f}^{\n}(\eta,\xi) + \frac{\psi(\eta,\xi)}{\norm{\g'(\xi)}}\hat{f}^{\ttau}_{\xi}(\eta,\xi),
   ~~\forall \x \in N_\Gamma(\epsilon).
    \label{eqn:div_f}
\end{eqnarray}
Lastly but very importantly, we replace $\mathbf{f}$ in \eqref{eqn:div_f} with $\nabla u$ and use \eqref{eqn:grad_u} to obtain
\begin{subequations}
    \begin{eqnarray}
\Delta u(\x) &=& \nabla\cdot (\nabla u(\x)) = \L(\hu(\eta, \xi)) ~~\text{with}\nonumber \\
\L(\hu(\eta, \xi)) &:=& \hu_{\eta\eta}(\eta,\xi)+J_0(\eta,\xi) \hu_{\xi\xi}(\eta,\xi) + J_1(\eta,\xi)\hu_{\eta}(\eta,\xi)  +J_2(\eta,\xi)\hu_{\xi}(\eta,\xi), \label{eqn:laplacian_formula}
\end{eqnarray}
where
\begin{equation}
\begin{split}
J_0(\eta,\xi) &=\left(\frac{\psi(\eta,\xi)}{\norm{\g'(\xi)}}\right)^2,\qquad J_1(\eta,\xi) =\kappa(\xi)\psi(\eta,\xi), \\
J_2(\eta,\xi) &=-\left(\frac{\psi(\eta,\xi)}{\norm{\g'(\xi)}}\right)^2\left(\eta \kappa'(\xi)\psi(\eta,\xi)+\frac{\g'(\xi)\cdot \g''(\xi)}{\norm{\g'(\xi)}^{2}}\right),
\end{split} \label{eqn:def_of_Js}
\end{equation}
and we will tacitly assume that the parametrization $\g(\xi)$ of the interface curve $\Gamma$ is $C^3$ whenever the term $\L(\hu(\eta, \xi))$ is involved.
\end{subequations}

{It is worth mentioning here that the interface $\Gamma$ might be described by a cartesian equation $F(\x)=0$ instead of a parametrization. In this case, we can obtain a parametrization $\g$ efficiently and accurately by solving the differential equation $\g'(\xi)=c(\xi)Q\nabla F(\g(\xi))$ where $Q$ is defined in \eqref{eqn:def_of_Q} and $c$ is a scalar function. For a detailed discussion of this topic and the various choices of $c$, we refer the reader to \cite{Kuznetsov2006,Kuznetsov2003}.

}


\commentout{
    and smooth enough: It is sufficient for $\Gamma$ to be $\color{red}C^2$ to construct the first order immersed discontinuous finite element space. On the other hand, we will require $\Gamma$ to be $C^3$ in order to construct higher order IDFE spaces, this requirement is enough to construct IDFE spaces of any order. However, {\color{red} in our analysis of the approximation capabilities, we will require $\Gamma$ to be smoother.}
}

\commentout{
        From our smoothness and connectedness assumptions, it follows that $\Gamma$ is locally diffeomorphic to an open interval and globally diffeomorphic to an interval or the unit circle \cite{milnor_1965}. However, a such global parametriztion might not be available or cheap. Instead, we will assume that for every interface element $K$, there is an interval $I_K=(a_K,b_K)$ and a regular parametrization $\g:I_K\to \Gamma$ such that $\Gamma \cap K\subset \tGamma_K:= \g(I_K)$. Here, the regularity of $\g$ means that $\g'(\xi)\ne 0$ for all $\xi\in I_K$.For instance, see \autoref{fig:rectangular_interface_element} for an example a such parametrization $\g$.
        \begin{figure}[ht]
            \centering
            \begin{tikzpicture}[scale=2]
        \draw [black]  (0,0) -- (0,1) -- (1.5,1) -- (1.5,0)   -- cycle;
          \draw[red] (.9,-.1) .. controls (.4,.4) and (.5,.5) .. (-.1,.6) node[midway,right]{$\tGamma_K$}
          ;
          \node at (0.4,.1) {$K^-$};
          \node at (.7,.8) {$K^+$};
          \fill[blue](0,0) circle[radius=.6pt];
          \fill[blue](1.5,0) circle[radius=.6pt];
          \fill[blue](1.5,1) circle[radius=.6pt];
           \fill[blue](0,1) circle[radius=.6pt];
        \end{tikzpicture}
        \caption{An example of an interface element $K$}\label{fig:rectangular_interface_element}
        \end{figure}

        At this point, we note that part of $\tGamma_K$ might lay outside of the element $K$. This is not a problem if $K$ is an interior element. However, if $K$ is a boundary element, then a portion $\tGamma_K$ may be outside of $\O$. Later in the numerical examples, we see that a such situation is not problematic, provided that we can compute $\g$.

        To simplify our calculation, we set $\ttau_K=|\g'|^{-1}\g'_K$ be the unit tangent vector where $|\cdot|$ is the Euclidean norm in $\mathbb{R}^2$ and we choose the unit normal vector to be

        Let $k= |\g'|^{-3} \left(\g'\btimes \g''\right)$, where $\u\btimes \v= \u^T Q\v$ is the 2D cross product, be the signed curvature. We note here that we are dropping the dependence on $K$ since $k$ is independent of the choice of the parametrization $\g$ (see \cite{do_carmo_differential_1976}). Lastly, we assume that the curvature $k$ is bounded in absolute value by $\kappa_\Gamma$.
}

\section{An IFE Space by Frenet Coordinates}  \label{sec:Frenet_space}

A key component of an IFE method is the construction of IFE functions on interface elements. While the construction of IFE functions based on linear or bilinear polynomials is
moderately involved \cite{guo_group_2019,he_approximation_2008,li_immersed_2004, li_new_2003,lin_rectangular_2001}, the construction of higher degree IFE functions gives rise to many challenges because of the impossibility for polynomials to satisfy the jump conditions across a general interface curve $\Gamma$. Some attempts for constructing IFE functions with higher degree polynomials are presented in \cite{benromdhane2014,adjerid_higher_2018,adjerid_high_2017,guo_higher_2019}. {Also, it is worth mentionning that the tubular neighborhood defined in \eqref{eq:tubular} was used in IFE methods in \cite{guo_immersed_2020,guzman_2017} in two and three dimensions}, but the information about the interface curve is insubstantially employed in either the line integral along the interface or the boundary of the subelements in weak forms for enforcing jump conditions.

As described in \autoref{sec:Frenet_coordinates}, for each interface element $K \in \mathcal{T}_h$, the Frenet transformation $P_\Gamma$ defined in \eqref{eq:paralle_curve} associates the fictitious element $K_F$ with its Frenet fictitious element $\hK_F$ in which the interface curve $\Gamma_{K_F}$ in $K_F$ is mapped to $\hGamma_{K_F} = R_\Gamma(\Gamma_{K_F})$ that is a line segment along the $\xi$ axis, see \autoref{fig:intf_fic_elements} for an illustration. This distinct feature opens the possibility to construct IFE functions with higher degree polynomials in the Frenet fictitious element $\hK_F$ because the transformed interface $\hGamma_{K_F}$ there is a simple line segment, and polynomials on both sides of $\hGamma_{K_F}$ can be easily constructed to form a piecewise polynomial that can satisfy the jump conditions across $\hGamma_{K_F}$, {\em precisely}.
The polynomial IFE functions are then transformed to the interface element $K$ by the Frenet transformation $P_\Gamma$ to form IFE functions on $K$. Of course, IFE functions
constructed in this way are similar to the isoparametric finite element functions in the sense that they are not necessarily polynomials; however, as to be shown later, these IFE functions will have desirable features such as  that they satisfy jump conditions \eqref{eqn:continuity_cond} and \eqref{eqn:normal_cond} {\em precisely} on the interface curve $\Gamma_K = \Gamma \cap K$ inside the interface element $K$.

\subsection{A local IFE space on every interface element}
Let $m\ge 1$ be an integer, and let $u$ be a function that is smooth enough and satisfies the jump conditions in \eqref{eqn:continuity_cond}, \eqref{eqn:normal_cond} and \eqref{eqn:laplacian_cond}. For an interface element $K\in \mathcal{T}_h$, we consider $\hu =u\circ P_\Gamma$ and let $\hbeta=\beta\circ P_\Gamma$. Then, by \eqref{eqn:normal_to_eta} and \eqref{eqn:laplacian_formula}, we can see that $\hu(\eta, \xi)$ satisfies the following jump conditions:
\begin{equation}
\bb{\hu}_{\hGamma_{K_F}}=0,
\qquad  \bb{\hbeta \hu_\eta }_{\hGamma_{K_F}}=0,
\qquad \bb{\hbeta
\frac{\p^j}{\p\eta^j}\L (\hu)}_{\hGamma_{K_F}}=0, ~~j=0,1,\dots,m-2.\label{eqn:flattened_interface_conditions}
\end{equation}
 We introduce two subelements of the Frenet fictitious
element $\hK_F$:
\begin{eqnarray}
\hK_F^- = \{(\eta, \xi) \in \hK_F, ~\eta <0\}, \quad \hK_F^+ = \{(\eta, \xi) \in \hK_F, ~\eta \geq 0\}
\label{eq:hKF_minus_plus}
\end{eqnarray}
which form a natural partition of $\hK_F$ according to the interface $\hGamma_{K_F}$ which, in turn,  is associated with the interface $\Gamma_{K_F}$ by the Frenet transformation. Without loss of generality, we will assume in the rest of the paper that $P_{\Gamma}(\hK^{\pm})\subset \Omega^{\pm}$.

Given a set $O \subset \mathbb{R}^2$, we let $Q^m(O)$ be the space of tensor product polynomials on $O$ of degree up to $m$ in each variable. Similarly, given a set $I\subset \mathbb{R}$, we let $\P^m(I)$ be the space of polynomials of degree up to $m$ on $I$. We then seek to construct IFE functions in the following piecewise polynomial form:
\begin{subequations}

\begin{equation}
\hphi(\eta,\xi)= \begin{cases}
    \displaystyle \sum_{0\le i\le m} \eta^i p^-_i(\xi),& \text{for~~} (\eta, \xi) \in \hK_F^-, ~p^-_i\in \P^m(\hGamma_{\hK_F}),\\ \\
    \displaystyle \sum_{0\le i\le m} \eta^i p^+_i(\xi),& \text{for~~} (\eta, \xi) \in \hK_F^+, ~p^+_i\in \P^m(\hGamma_{\hK_F}),
\end{cases}
\end{equation}
that can satisfy the interface jump condition \eqref{eqn:flattened_interface_conditions}. We can rewrite
$\hphi(\eta,\xi)$ as follows:
\begin{eqnarray}
\hphi(\eta, \xi) &=& \hphi^-(\eta, \xi) \chi_{\hK_F^-}(\eta, \xi) + \hphi^+(\eta, \xi) \chi_{\hK_F^+}(\eta, \xi), \label{eqn:phi_p_m_decomp} \\[5pt]
\hphi^s(\eta, \xi) &=& \sum_{0\le i\le m} \eta^i p^s_i(\xi), ~~p^s_i\in \P^m(\hGamma_{K_F}), ~s = -, +, \label{eqn:phi_p_m_decomp_1}
\end{eqnarray}
\end{subequations}
where $\chi_{\hK_F^s}(\eta, \xi)$ is the characteristic function for the set $\hK_F^s, ~s = -, +$.

While the application of the conditions $\bb{\hphi}_{\hGamma_{K_F}} = 0$ and $\bb{\hbeta \hphi_{\eta}}_{\hGamma_{K_F}}=0$ leads to useful algebraic equations for determining $\hphi^s(\eta, \xi), ~s = -, +$, we cannot enforce the last jump condition in \eqref{eqn:flattened_interface_conditions} exactly unless $\Gamma_{K_F}$ is an arc with a {constant} curvature, such as a circle or a line. In general, the curvature $\kappa(\xi)$ of $\Gamma$ may be a non-polynomial function and so is $\L(\hphi)$. Hence, the last jump condition in \eqref{eqn:flattened_interface_conditions} will not provide polynomial equations for us to directly determine the polynomial terms $p^\pm_i(\xi), i = 2, 3, \ldots, m$ in $\hphi(\eta, \xi)$. For this reason we impose this jump condition on $\hphi$ weakly and use the following jump conditions whenever necessary:
\begin{subequations}\label{eqn:extended_laplacian_IFE_cond}
\begin{eqnarray}
\int_{\hGamma_{K_F}}\bb{\hphi}_{\hGamma_{K_F}} v =0, \qquad  \int_{\hGamma_{K_F}}\bb{\hbeta \hphi_\eta }_{\hGamma_{K_F}}v =0, ~~\forall v \in \P^m(\hGamma_{K_F}), \label{eqn:extended_laplacian_IFE_cond_a}
\end{eqnarray}
\begin{equation}
\int_{\hGamma_{K_F}} \bb{\frac{\p^j}{\p\eta^j}\L(\hphi)}_{\hGamma_{K_F}} v =0, \qquad \forall v \in \P^m(\hGamma_{K_F}),\qquad j=0,1,\dots,m-2. \label{eqn:extended_laplacian_IFE_cond_b}
\end{equation}
\end{subequations}
We note that \eqref{eqn:extended_laplacian_IFE_cond_a} is equivalent to the first two jump conditions in \eqref{eqn:flattened_interface_conditions}.

The following lemma describes how the components $\hphi^s(\eta, \xi), ~s = -, +$ of the piecewise polynomial
$\hphi(\eta, \xi)$ specified by \eqref{eqn:phi_p_m_decomp} are related to each other by the interface jump conditions.

\begin{lemma} \label{lem:existence_uniqueness_extension}
Given a polynomial $\hphi^{\pm} \in  Q^m(\hK_F^\pm)$, there is a unique polynomial $\hphi^{\mp} \in  Q^m(\hK_F^\mp)$ such that
the piecewise $Q^m$ polynomial $\hphi = \hphi^- \chi_{\hK_F^-} + \hphi^+ \chi_{\hK_F^+}$ satisfies the jump conditions \eqref{eqn:extended_laplacian_IFE_cond}.
\end{lemma}
\begin{proof}
Without loss of generality, we assume that $\hphi^- \in Q^m(\hK_F^-)$ is given. By varying $v$ among a basis
of $\P^m(\hGamma_{K_F})$, we can see that jump conditions \eqref{eqn:extended_laplacian_IFE_cond} lead to a square system of linear equations about the coefficients of $\hphi^+$. Hence, we can assume $\hphi^- = 0$ and proceed to show that $\hphi^+ = 0$ provided that $\hphi$ in the format specified in \eqref{eqn:phi_p_m_decomp} satisfies the jump conditions \eqref{eqn:extended_laplacian_IFE_cond}. By \eqref{eqn:phi_p_m_decomp_1} and the fist jump condition in \eqref{eqn:extended_laplacian_IFE_cond_a}, we have
\begin{eqnarray*}
\int_{\hGamma_{K_F}}p_0^+(\xi)v d\xi = \int_{\hGamma_{K_F}}\bb{\hphi}_{\hGamma_{K_F}}v  = 0, ~~\forall v \in \P^m(\hGamma_{K_F})
\end{eqnarray*}
which means $p_0^+$ is the zero polynomial and
\begin{eqnarray*}
\hphi^+(\eta, \xi) &=& \sum_{1\le i\le m} \eta^i p^+_i(\xi).
\end{eqnarray*}
Then, by the second jump condition in \eqref{eqn:extended_laplacian_IFE_cond_a}, we have
\begin{eqnarray*}
\int_{\hGamma_{K_F}}\beta^+ p_1^+(\xi) v d\xi= \int_{\hGamma_{K_F}}\bb{\hbeta \hphi_{\eta}}_{\hGamma_{K_F}} v d\xi=0, ~~\forall v \in \P^m(\hGamma_{K_F})
\end{eqnarray*}
which implies $p_1^+$ is the zero polynomial. The proof is finished if $m = 1$. If $m \geq 2$, we note that
\begin{eqnarray*}
\hphi^+(\eta,\xi)= \sum_{i=2}^m\eta^i p_i^+(\xi),\qquad \text{which implies}~~ \hphi_{\xi}^+(0,\xi)=\hphi_{\xi\xi}^+(0,\xi)=\hphi^+_{\eta}(0,\xi) = 0.
\end{eqnarray*}
Then, by \eqref{eqn:laplacian_formula} and the jump condition \eqref{eqn:extended_laplacian_IFE_cond_b} with
$j = 0$, we have
\begin{eqnarray*}
\int_{\hGamma_{K_F}}2\beta^+p_2^+(\xi)v d\xi = \int_{\hGamma_{K_F}}\beta^+\hphi_{\eta \eta}^+(0,\xi) v d\xi = \int_{\hGamma_{K_F}}\bb{\hbeta\L (\hphi)}_{\hGamma_{K_F}}v d\xi =0,
~~\forall v \in \P^m(\hGamma_{K_F}),
\end{eqnarray*}
which means that $p_2^+$ is the zero polynomial. Next, we proceed by mathematical induction and assume  $p_0^+=p_1^+=\dots=p_{i_0}^+=0$ for some integer $2\le i_0< m$. Then, we have
\begin{eqnarray}
&&\hphi^+(\eta,\xi)=\sum_{i=i_0+1}^m \eta^i p_i^+(\xi), \label{eqn:induction_hyp} \\
\text{and }&& \frac{\partial^j}{\partial \eta^j}\hphi^+(0,\xi)=\frac{\partial^j}{\partial \eta^j}\hphi_{\xi}^+(0,\xi)=\frac{\partial^j}{\partial \eta^j}\hphi_{\xi\xi}^+(0,\xi)= 0,\qquad j=0,1,\dots,i_0.
\label{eqn:induction_hyp_1}
\end{eqnarray}
Also, by \eqref{eqn:laplacian_formula}, we have
\begin{align}
\begin{split}\label{eqn:L_big_sum}
&\frac{\partial^{i_0-1}}{\partial\eta^{i_0-1}}\L(\hphi^+) = \frac{\partial^{i_0+1}}{\partial\eta^{i_0+1}}\hphi^+
+\frac{\partial^{i_0-1}}{\partial\eta^{i_0-1}}\left(J_0\hphi^+_{\xi\xi}+ J_1\hphi^+_{\eta} +J_2 \hphi^+_{\xi}\right), \\
&=\frac{\partial^{i_0+1}\hphi^+}{\partial\eta^{i_0+1}}
+\sum_{l=0}^{i_0-1}\binom{i_0-1}{l}\left( \frac{\partial^{i_0-1-l}J_0}{\partial\eta^{i_0-1-l}}\frac{\partial^l \hphi^+_{\xi\xi}}{\partial\eta^l}+
 \frac{\partial^{i_0-1-l}J_1}{\partial\eta^{i_0-1-l}}\frac{\partial^{l+1}\hphi^+}{\partial\eta^{l+1}}
 +
 \frac{\partial^{i_0-1-l}J_2}{\partial\eta^{i_0-1-l}}\frac{\partial^l \hphi^+_{\xi}}{\partial\eta^l}\right).
\end{split}
\end{align}
Then, applying \eqref{eqn:induction_hyp}, \eqref{eqn:induction_hyp_1}, and \eqref{eqn:L_big_sum}
to the jump condition \eqref{eqn:extended_laplacian_IFE_cond_b} with $j = i_0-1$, we have
\begin{eqnarray*}
\int_{\hGamma_{K_F}}\beta^+(i_0+1)!p_{i_0+1}^+(\xi) v d\xi &=& \int_{\hGamma_{K_F}}\beta^+\frac{\partial^{i_0-1}}{\partial\eta^{i_0-1}}\L(\hphi^+)(0,\xi) v d\xi \\
&=& \int_{\hGamma_{K_F}} \bb{\hbeta
\frac{\p^{i_0-1}}{\p\eta^{i_0-1}}\L (\hphi)}_{\hGamma_{K_F}}v d\xi=0, ~~\forall v \in \P^m(\hGamma_{K_F}),
\end{eqnarray*}
which allows us to conclude that $p_{i_0+1}^+$ is the zero polynomial. By mathematical induction, we should have
$p_i^+, 0 \leq i \leq m$ are zero polynomials; hence, we conclude that $\hphi^+ = 0$.

\end{proof}


\autoref{lem:existence_uniqueness_extension} means there exist IFE functions based on $Q^m$ polynomials that satisfy the transformed jump conditions \eqref{eqn:extended_laplacian_IFE_cond} in the Frenet fictitious element $\hK_F$ because we can pick up any nonzero polynomials $\hphi^s \in Q^m(\hK_F^s)$, the lemma guarantees the existence of another polynomial $\hphi^{s'} \in Q^m(\hK_F^{s'})$, where $s' = \pm$ when $s = \mp$, such that $\hphi = \hphi^- \chi_{\hK_F^-} + \hphi^+ \chi_{\hK_F^+}$ is the desired IFE function satisfying the transformed jump condition \eqref{eqn:extended_laplacian_IFE_cond}. Consequently, the following IFE space on the Frenet fictitious element $\hK_F$ is  well-defined:
\begin{eqnarray}
\hV^m_{\hbeta}(\hK_F) &=&\left\{ \hphi:\hK_F\to \mathbb{R}\mid
\hphi|_{\hK_F^{\pm}} \in Q^m(\hK_F^{\pm})~~\text{ and } ~\hphi\text{ satisfies }\eqref{eqn:extended_laplacian_IFE_cond}\right\}. \label{eqn:def_IFE_space}
\end{eqnarray}

On the interface element $K\in \mathcal{T}_h^i$ associated with $\hK_F$, we can use the inverse Frenet transform
$R_\Gamma$ and the IFE space $\hV^m_{\hbeta}(\hK_F)$ on the Frenet fictitious element to define an IFE space based on $Q^m$ polynomials as follows:

\begin{equation}
\V^m_{\beta}(K)=\left\{\hphi\circ(R_\Gamma)|_{K} \mid \hphi \in \hV^m_{\hbeta}(\hK_F)\right\}.
\label{eqn:def_IIFE_physical}
\end{equation}
For a general interface  $\Gamma$, an IFE function $\phi \in \V^m_{\beta}(K)$ is not necessarily a
polynomial because of the nonlinearity of $R_\Gamma$ described in \eqref{eq:P_inv}. Nevertheless, IFE functions
defined by \eqref{eqn:def_IIFE_physical} are constructed with $Q^m$ polynomials and the Frenet transformation with fundamental differential geometry features of the interface curve $\Gamma$, and for a general interface curved
$\Gamma$, they possess the following distinct features that are not present in any of the IFE functions proposed in the literature.

\begin{lemma}\label{lem:conformity}
    Let $K$ be an interface  element. Then every IFE function $\phi\in \V^m_{\beta}(K)$ satisfies
    \begin{enumerate}
        \item the point-wise continuity across the interface curve $\Gamma_K = \Gamma \cap K$, {\it i.e.}, $\bb{\phi}_{\Gamma_K}=0$,
        \item the point-wise flux continuity across the interface curve $\Gamma_K = \Gamma \cap K$, {\it i.e.}, $\bb{\beta\frac{\p}{\p \n}\phi}_{\Gamma_K}=0$.
    \end{enumerate}
\end{lemma}

\begin{proof}
Let $\phi \in \V^m_{\beta}(K)$. Then there exists $\hphi \in \hV^m_{\hbeta}(\hK_F)$ such that
\begin{eqnarray*}
\phi(\x) = \hphi(R_\Gamma(\x)) = \hphi(\eta, \xi).
\end{eqnarray*}
Then, Claim 1. follows from the fact that both $\hphi$ and $R_\Gamma(\x)$ are continuous functions. Claim 2. follows
from the facts that $\frac{\p}{\p \n}\phi$ satisfies \eqref{eqn:normal_to_eta} and
$\hphi$ satisfies the jump condition across $\hGamma_{K}=\hGamma_{K_F}\cap \hK$ so that we have
\begin{eqnarray*}
\bb{\beta\frac{\p}{\p \n}\phi(\x)}_{\Gamma_K} = \bb{\hbeta \hphi_{\eta}(\eta, \xi)}_{(\eta, \xi) \in \hGamma_{K_F}} = 0.
\end{eqnarray*}
We used the fact that $K \subseteq K_F$ and $R_\Gamma(\Gamma_K) \subseteq \hGamma_{K_F}$ consequently.

\end{proof}

\subsection{A basis for the local IFE space}
\label{subsec:basis_IFE}

While the proof of \autoref{lem:existence_uniqueness_extension} provides some hints on how to construct the IFE shape functions, we now present a specific procedure for generating a basis for the local IFE space $\V^m_{\beta}(K)$ on every
interface element $K \in \mathcal{T}_h^i$ which can be used in computations involving the proposed IFE spaces.
By the definition of $\V^m_{\beta}(K)$, we only need to construct a basis for the $Q^m$ IFE space $\hV^m_{\hbeta}(\hK_F)$
on the Frenet fictitious element $\hK_F$. By \autoref{lem:existence_uniqueness_extension}, we know that
$\dim(\hV^m_{\hbeta}(\hK_F)) = (m+1)^2$ which is actually the dimension of $Q^m(\hK_F)$.

As another major advantage of our proposed approach to use the Frenet apparatus of the interface curve $\Gamma$ in IFE functions, we can easily construct a majority of the basis functions for $\hV^m_{\hbeta}(\hK_F)$ without any computations. In fact, the first $m(m+1)$ basis functions are available almost freely as stated in \autoref{lem:partial_basis} below. Let $\{p_i(\xi)\}_{i = 0}^m$ be a basis of $\P^m(\hGamma_{K_F})$, then it can be easily verified that the following polynomials form a basis of $Q^m(\hK_F)$:
\begin{eqnarray}
\begin{matrix}
p_0(\xi), &\eta p_0(\xi), &\N_1(\eta, \xi) = \eta^2 p_0(\xi), &\cdots &\N_{m^2-m-1}(\eta, \xi) = \eta^m p_0(\xi) \\
p_1(\xi), &\eta p_1(\xi), &\N_2(\eta, \xi) = \eta^2 p_1(\xi), &\cdots &\N_{m^2-m}(\eta, \xi) = \eta^m p_1(\xi) \\
p_2(\xi), &\eta p_2(\xi), &\N_3(\eta, \xi) = \eta^2 p_2(\xi), &\cdots &\N_{m^2-m+1}(\eta, \xi) = \eta^m p_2(\xi) \\
\vdots    &\vdots         &\vdots          &\cdots &\vdots \\
p_m(\xi),&\eta p_m(\xi),  &\N_m(\eta, \xi) = \eta^2 p_m(\xi), &\cdots & \N_{m^2-1}(\eta, \xi) = \eta^m p_m(\xi)
\end{matrix} \label{eq:Q^m_basis_general}
\end{eqnarray}

\begin{lemma}\label{lem:partial_basis}
    Let $K$ be an interface element, then for all integers $0\le i\le m$ and $1\le j\le m$, the functions \begin{equation}
    \hphi_{i,j}(\eta,\xi) := \frac{1}{\hbeta(\eta,\xi)}\eta^jp_i(\xi), ~~0 \leq i \leq m, ~1 \leq j \leq m, \label{eqn:hphi_i_j_lem}
    \end{equation}
     are in the IFE space $\hV^m_{\beta}(\hK_F)$. Furthermore, the functions $\{\hphi_{i,j}\}_{i=0,j=1}^m$ are linearly independent.
\end{lemma}

\begin{proof}
For every $\hphi = \hphi_{i,j}$ defined by \eqref{eqn:hphi_i_j_lem}, it is easily to see that
$\hphi|_{\hK_F^s} \in Q^m(\hK_F^s)$ and $\hphi(0^+,\xi)=0=\hphi(0^-,\xi)$ since $j\ge 1$; consequently, we have $\bb{\hphi}_{\hGamma_{K_F}}=0$ because $\hGamma_{K_F}$ is a subset of the $\xi$ axis. Next, by definition, $\hbeta\hphi$ is a polynomial on $\hK_F$; hence, we have $\bb{\hbeta \hphi_{\eta}}_{\hGamma_{K_F}}=0$. Moreover, by formula \eqref{eqn:laplacian_formula} and the fact that $\hbeta\hphi$ is polynomial, we know that
$\L (\hbeta \hphi)$ is a rational function of $\eta$ so that $\frac{\p^j}{\p\eta^j}\L (\hbeta \hphi)$ is continuous on the Frenet fictitious element $\hK_F$. Therefore, by the linearity of $\L$, we have
\begin{eqnarray*}
\bb{\hbeta \frac{\p^l}{\p\eta^l}\L (\hphi)}_{\hGamma_{K_F}} =
\bb{\frac{\p^l}{\p\eta^l}\L (\hbeta \hphi)}_{\hGamma_{K_F}} = 0, ~~l=0,1,\dots,m-2.
\end{eqnarray*}
In summary, we have shown that $\hphi_{i, j}$ is in the IFE space $\hV^m_{\beta}(\hK_F)$, for all $0\le i\le m $ and $1\le j\le m$.

Lastly, the set $\{\hphi_{i,j}^+\}_{i=0,j=1}^m$ is linearly independent in $Q^m(\hK_F^+)$. Hence, the set of functions $\{\hphi_{i,j}\}_{i=0,j=1}^m$ is linearly independent.

\end{proof}


\commentout{
        Furthermore, we can see that $\hphi_{i,j}$ as defined in \autoref{lem:partial_basis} satisfies the interface conditions \eqref{eqn:flattened_interface_conditions} exactly. In other words, the functions $\phi_{i,j}=\hphi\circR_\Gamma$ behave similarly to the solution $u$ at the interface since they  the jump conditions \eqref{eqn:continuity_cond}, \eqref{eqn:normal_cond} and $\eqref{eqn:laplacian_cond}$.
}

We then proceed to the construction of the remaining $m+1$ functions for a basis of the IFE space $\hV^m_{\beta}(\hK_F)$.

\begin{lemma} \label{lem:low_order_extension}
    Let $i$ be an integer such that $0\le i\le m$ and let $\phi^-_{i,0}(\eta,\xi)=p_i(\xi)$ on $\hK^-$, then there is a unique $\phi_{i,0}^+\in Q^m(\hK^+)$ such that $\hphi_{i,0}\in \hV^m_{\beta}(\hK)$.
\end{lemma}

\begin{proof}
    This is a direct consequence of \autoref{lem:existence_uniqueness_extension}.
\end{proof}
By combining the previous two lemmas, we obtain a complete basis of $\hV^m_{\beta}(\hK)$ as stated in the following corollary.

\begin{corollary} \label{cor:basis_ref}
    Let  $B=\{\hphi_{i,j}\}_{i,j=0}^m$ be the set of functions defined in Lemmas \ref{lem:partial_basis} and \ref{lem:low_order_extension}, then $B$ is a basis of $\hV^m_{\beta}(\hK_F)$.
\end{corollary}


Even though \autoref{lem:low_order_extension} (or \autoref{lem:existence_uniqueness_extension}) is related to the existence of the first $m+1$ functions $\hphi_{i, 0}(\eta, \xi), ~0 \leq i \leq m$ in a basis for the IFE space $\hV^m_{\beta}(\hK_F)$, it actually shed light on how to construct these functions, and we now present the associated procedure. Our task is to start from $\hphi_{i, 0}^-(\eta,\xi)= p_i(\xi), ~0 \leq i \leq m$ and we construct $\hphi_{i, 0}^+\in Q^m(\hK)$ such that $\hphi_{i, 0} = \hphi_{i, 0}^- \chi_{\hK_F^-} + \hphi_{i, 0}^+ \chi_{\hK_F^+} \in \hV^m_{\beta}(\hK_F)$.

First, following  a similar argument to the one used in the proof of \autoref{lem:existence_uniqueness_extension}, we have $\hphi_{i, 0}^+(0,\xi)=p_i(\xi)$ and $\frac{\p}{\p\eta} \hphi_{i, 0}^+(0,\xi)=0$. Therefore, $\hphi_{i, 0}^+$ can be written in terms of the basis functions listed in \eqref{eq:Q^m_basis_general} as follows:

\begin{equation}
\hphi_{i, 0}^+(\eta,\xi)= p_i(\xi) + \sum_{l=1}^{m^2-1}c_l^{(i)}\N_l(\eta,\xi), \label{eqn:phip_j_2_single_sum}
\end{equation}
with its coefficients $\c^{(i)} = (c_l^{(i)})_{l = 1}^{m^2-1}$ to be determined. Next, we use the weak jump conditions stated in \eqref{eqn:extended_laplacian_IFE_cond} with the test function $v = p_k(\xi), k = 0, 1, \dots, m$ to have

\begin{eqnarray}
&&\sum_{l=1}^{m^2-1}c_l^{(i)}\int_{\hGamma_{K_F}}\frac{\partial^j}{\partial \eta^j} \L(\N_l)(0,\xi) p_k(\xi)d\xi
=\frac{\beta^--\beta^+}{\beta^+} \int_{\hGamma_{K_F}} \frac{\partial^j}{\partial \eta^j} \L(p_i(\xi))(0,\xi)p_k(\xi) d\xi, \label{eqn:expanded_ife_system}\\
&&\hspace{1.5in}\text{for}\hspace{0.1in} 0\le k\le m, ~~0\le j\le m-2. \nonumber
\end{eqnarray}

The equations in \eqref{eqn:expanded_ife_system} lead to the following linear system for computing the coefficients of
$\hphi_{i, 0}^+(\eta, \xi)$ described in \eqref{eqn:phip_j_2_single_sum}:
\begin{equation}
\mathbf{A}\c^{(i)} =\frac{\beta^- - \beta^+}{\beta^+} \b(i) \label{eqn:compact_ife_system}
 \end{equation}
where
\begin{eqnarray}
\begin{aligned} \label{eqn:def_A}
&\mathbf{A}=\begin{bmatrix}
    A^{(0)} \\ A^{(1)} \\ \vdots\\ A^{(m-2)}
\end{bmatrix},\qquad \b(i) = \begin{bmatrix}
\b^{(0)}(i) \\
\b^{(1)}(i) \\
\vdots \\
\b^{(m-2)}(i)
\end{bmatrix}, \\
&A_{k, l}^{(j)}=\int_{\hGamma_{K_F}}\frac{\partial^j}{\partial \eta^j} \L(\N_l)(0,\xi)p_k(\xi) d\xi,\qquad \b^{(j)}_k(i)=   \int_{\hGamma_{K_F}} \frac{\partial^j}{\partial \eta^j} \L(p_i(\xi))(0,\xi)p_k(\xi) d\xi, \\
& 0 \leq j \leq m-2, ~~1 \leq l \leq m^2-1, ~~0 \leq k \leq m.
\end{aligned}
\end{eqnarray}
We note that the matrix $\mathbf{A}$ does not depend on $\beta^+$ or $\beta^-$ and its invertibility is guaranteed by \autoref{lem:low_order_extension}.

In summary, for the construction of $\hphi_{i, 0}(\eta, \xi)$, we chose $\hphi_{i, 0}^-(\eta, \xi) = p_i(\xi)$, solve the linear system \eqref{eqn:compact_ife_system} for $\c^{(i)}$, form $\hphi_{i, 0}^+(\eta, \xi)$ by \eqref{eqn:phip_j_2_single_sum}, then
set $\hphi_{i, 0} = \hphi_{i, 0}^- \chi_{\hK_F^-} + \hphi_{i, 0}^+ \chi_{\hK_F^+} \in \hV^m_{\beta}(\hK_F)$.

By  \autoref{cor:basis_ref}, we can express the IFE space associated with an interface element
$K \in \mathcal{T}_h^i$ as the span of the basis functions constructed above in the following ways:
\begin{eqnarray}
\begin{aligned}\label{eq:IFE_space_span}
\hV^m_{\beta}(\hK_F) &= \operatorname{span}\{\hphi_{i, j}(\eta, \xi), ~0 \leq i, j \leq m\}, \\
\V^m_{\beta}(K) &= \operatorname{span}\{\phi_{i,j}(\x) = \hphi_{i, j}\circ R_\Gamma(\x), ~0 \leq i, j \leq m\}.
\end{aligned}
\end{eqnarray}

\noindent
{\bf Remarks}: Here are a few features of the proposed procedure for constructing a basis for the local IFE space on an interface element $K \in \mathcal{T}_h^i$:

\begin{itemize}
    \item
    By \autoref{lem:conformity}, all the basis functions $\phi_{i, j}(\x), ~~0 \leq i, j \leq m$ of the IFE space $\V^m_{\beta}(K)$ satisfy the jump conditions  \eqref{eqn:continuity_cond} and \eqref{eqn:normal_cond} for the interface problem.
    Moreover, by \autoref{lem:partial_basis}, when $m \geq 2$, the majority of the basis functions also satisfy the extended jump conditions specified in \eqref{eqn:laplacian_cond}, these basis functions are
    $\phi_{i, j}(\x), ~~0 \leq i \leq m, 1 \leq j \leq m$. This is a notable feature that are not present in any of the IFE functions proposed in the literature.

     \item
     If $\beta^+=\beta^-$, then $\c^{(i)}=\mathbf{0}$. Consequently, we have $\hV^m_{\beta}(\hK_F)=Q^m(\hK_F)$ which implies the consistency of the proposed IFE space with the standard $Q^m$ finite element space.

    \item
    We can replace the basis $\{p_i(\xi)\}_{i = 0}^m$ of $\P^m(\hGamma_{K_F})$ by any other basis with desirable features. For example, we can use the basis $\{\xi^i\}_{i = 0}^m$ for its popularity and simplicity. Or, for orthogonality, we can use $\{L_i(\xi)\}_{i=0}^m$, where $L_i(\xi)$ is the $i$-th degree Legendre polynomial on $\hGamma_{K_F}$.

\end{itemize}

\subsection{Computational aspects of the construction of an IFE basis}\label{subsec:comp_aspects_construction}

In this subsection, we discuss a few computational aspects about constructing an IFE basis on an interface element
$K\in \mathcal{T}_h^i$. By \autoref{lem:partial_basis}, the construction for most of the basis functions for the IFE space $\V^m_{\beta}(K)$ costs almost nothing. Moreover, the construction cost for $\phi_{i, 0}(\x), ~0 \leq i \leq m$ is moderate. This is because the matrix $\mathbf{A}$ in \eqref{eqn:compact_ife_system} is the same for every $0\le i\le m$ which can be used in the construction of all $\hphi_{i, 0}(\eta, \xi), ~~0\leq i \leq m$. This is advantageous for constructing higher degree IFE functions. In addition, with $\{\N_l\}_{l=1}^{m^2-1}$ specified in \eqref{eq:Q^m_basis_general}, the matrix $\mathbf{A}$ will be block lower triangular, a feature that can be used to reduce the cost when solving \eqref{eqn:compact_ife_system} for the coefficient vector $\c^{(i)}$.

For setting up the matrix $\mathbf{A}$ and the vector $\b(i)$ in the linear system \eqref{eqn:compact_ife_system}, we need to evaluate the following seemingly complex formula:
\begin{eqnarray}
\begin{aligned} \label{eq:eta_d_L}
&\frac{\partial^j}{\partial \eta^j} \L(v)(0,\xi) = \frac{\partial^{j+2}v(0,\xi)}{\partial\eta^{j+2}} \\
& +\sum_{l=0}^{j}\binom{j}{l}\left(
\frac{\partial^{l}J_0(0,\xi)}{\partial\eta^{l}}
\frac{\partial^{j-l} v_{\xi\xi}(0,\xi)}{\partial\eta^{j-l}}+
 \frac{\partial^{l}J_1(0,\xi)}{\partial\eta^{l}}
 \frac{\partial^{j-l+1}v(0,\xi)}{\partial\eta^{j-l+1}}+
 \frac{\partial^{l}J_2(0,\xi)}{\partial\eta^{l}}
 \frac{\partial^{j-l} v_{\xi}(0,\xi)}{\partial\eta^{j-l}}\right)
\end{aligned}
\end{eqnarray}
with $v = \N_l$ or $v = p_i$. However, all the terms in this formula can be easily prepared. First, since $v$ is a polynomial, its derivatives
\begin{eqnarray*}
\frac{\partial^{j+2}v(0,\xi)}{\partial\eta^{j+2}}, ~\frac{\partial^{j-l} v_{\xi\xi}(0,\xi)}{\partial\eta^{j-l}}, ~\frac{\partial^{j-l+1}v(0,\xi)}{\partial\eta^{j-l+1}},
~\frac{\partial^{j-l} v_\xi(0,\xi)}{\partial\eta^{j-l}}
\end{eqnarray*}
can be obtained easily. For the $\eta$ partial derivatives of $J_0, J_1$ and $J_2$, first, by \eqref{eqn:def_of_Js}, we have
\begin{eqnarray*}
&&\frac{\partial^l J_0(\eta, \xi)}{\partial \eta^l} = \frac{1}{\norm{\g'(\xi)}^2}\frac{\partial^l \psi^2(\eta, \xi)}{\partial \eta^l}, ~~\frac{\partial^l J_1(\eta, \xi)}{\partial \eta^l} = \kappa(\xi) \frac{\partial^l \psi(\eta, \xi)}{\partial \eta^l} \\
&&\frac{\partial^l J_2(\eta, \xi)}{\partial \eta^l} = -\frac{\kappa'(\xi)}{\norm{\g'(\xi)}^2}\left(
\eta \frac{\partial^l \psi^3(\eta, \xi)}{\partial \eta^l} + l \frac{\partial^{l-1} \psi^3(\eta, \xi)}{\partial \eta^{l-1}} \right) \\
&& \hspace{1.5in} - \frac{\g'(\xi)\cdot \g''(\xi)}{\norm{\g'(\xi)}^4} \frac{\partial^l \psi^2(\eta, \xi)}{\partial \eta^l}.
\end{eqnarray*}
Then, by the following expansions for the powers of $\psi(\eta,\xi)=(1+\kappa(\xi)\eta)^{-1}$:
\begin{eqnarray*}
&&\psi(\eta,\xi) =\sum_{j=0}^{\infty}(-1)^j\kappa(\xi)^j\eta^j,\qquad
\left(\psi(\eta,\xi)\right)^2 =\sum_{j=0}^{\infty}(-1)^j (j+1)\kappa(\xi)^j\eta^j, \\
&&\left(\psi(\eta,\xi)\right)^3 =\sum_{j=0}^{\infty}(-1)^j \frac{(j+1)(j+2)}{2}\kappa(\xi)^j\eta^j,
\end{eqnarray*}
we have:
\begin{eqnarray*}
&&\frac{\p^l\psi(0,\xi)}{\p\eta^l}= (-1)^l l!\kappa(\xi)^l,\quad
\frac{\p^l\psi^2(0,\xi)}{\p\eta^l}= (-1)^l (l+1)!\kappa(\xi)^l, \\
&&\frac{\p^{l-1}\psi^3(0,\xi)}{\p\eta^{l-1}}= (-1)^{l-1} \frac{(l+1)!}{2}\kappa(\xi)^{l-1}, \quad \frac{\p^l\psi^3(0,\xi)}{\p\eta^l}= (-1)^l \frac{(l+2)!}{2}\kappa(\xi)^l.
\end{eqnarray*}
Therefore, we can evaluate the $\eta$ partial derivatives of $J_0, J_1$ and $J_2$ on the $\xi$ axis by the following much simplified formulas:

\begin{eqnarray}
\begin{aligned}\label{eq:J_derivatives}
&\frac{\partial^l J_0(0, \xi)}{\partial \eta^l} = \frac{(-1)^l (l+1)!\kappa(\xi)^l}{\norm{\g'(\xi)}^2}, \quad \frac{\partial^l J_1(0, \xi)}{\partial \eta^l} = \kappa(\xi)(-1)^l l!\kappa(\xi)^l, \\
&\frac{\partial^l J_2(0, \xi)}{\partial \eta^l} = -\frac{\kappa'(\xi)}{\norm{\g'(\xi)}^2}\left(
l (-1)^{l-1} \frac{(l+1)!}{2}\kappa(\xi)^{l-1} \right) - \frac{\g'(\xi)\cdot \g''(\xi)}{\norm{\g'(\xi)}^4}(-1)^l (l+1)!\kappa(\xi)^l.
\end{aligned}
\end{eqnarray}
Hence, the evaluation of formula \eqref{eq:eta_d_L} is actually straightforward because all the terms in \eqref{eq:eta_d_L} have explicit formulas. The simplicity of these formulas is another benefit of using the Frenet apparatus of the interface curve $\Gamma$.

In the literature for higher degree IFE methods, it has been noted that the linear system for determining the coefficients of an IFE basis function on an interface element $K\in \mathcal{T}_h^i$ can be severely ill-conditioned  when
$\min(|\hK^+|,|\hK^-|)\ll \max(|\hK^+|,|\hK^-|)$, that is, when
$K$ has an extremely small
sub-element formed by the interface $\Gamma\cap K$. A remedy to alleviate this issue is to construct the IFE basis functions
on a fictitious element \cite{adjerid_high_2017,guo_higher_2019}. However the fictitious element proposed in the literature requires a user-chosen parameter. This parameter acts like most  parameters for regularization whose best choice is often uncertain.

The construction method for IFE basis functions on an interface element $K\in \mathcal{T}_h^i$ proposed in this article also depends on a fictitious element $K_F$. However, the fictitious element $K_F$ here acts in a self-regulatory way in the sense that $K_F$ is parameter free. Moreover, the matrix $\mathbf{A}$ and the vector $\b(i)$ in the linear system \eqref{eqn:compact_ife_system} for determining the coefficients of the IFE basis functions $\hphi_{i, 0}(\eta, \xi), ~0\leq i \leq m$ are constructed on the Frenet fictitious interface
$\hGamma_{K_F}$ predetermined by the element $K$ itself and is not susceptible to the small sub-element cut from the $K$ by the interface curve $\Gamma \cap K$, implying that the condition of the linear system \eqref{eqn:compact_ife_system} for determining the coefficients of the IFE basis functions will not be negatively affected by the small sub-element issue. Hence, the proposed construction for the IFE basis functions should be robust, and this feature is  another advantage of employing the Frenet apparatus of the interface curve $\Gamma$ in IFE methods.

\begin{figure}[htbp]
    \centering
    \includegraphics[scale=.3]{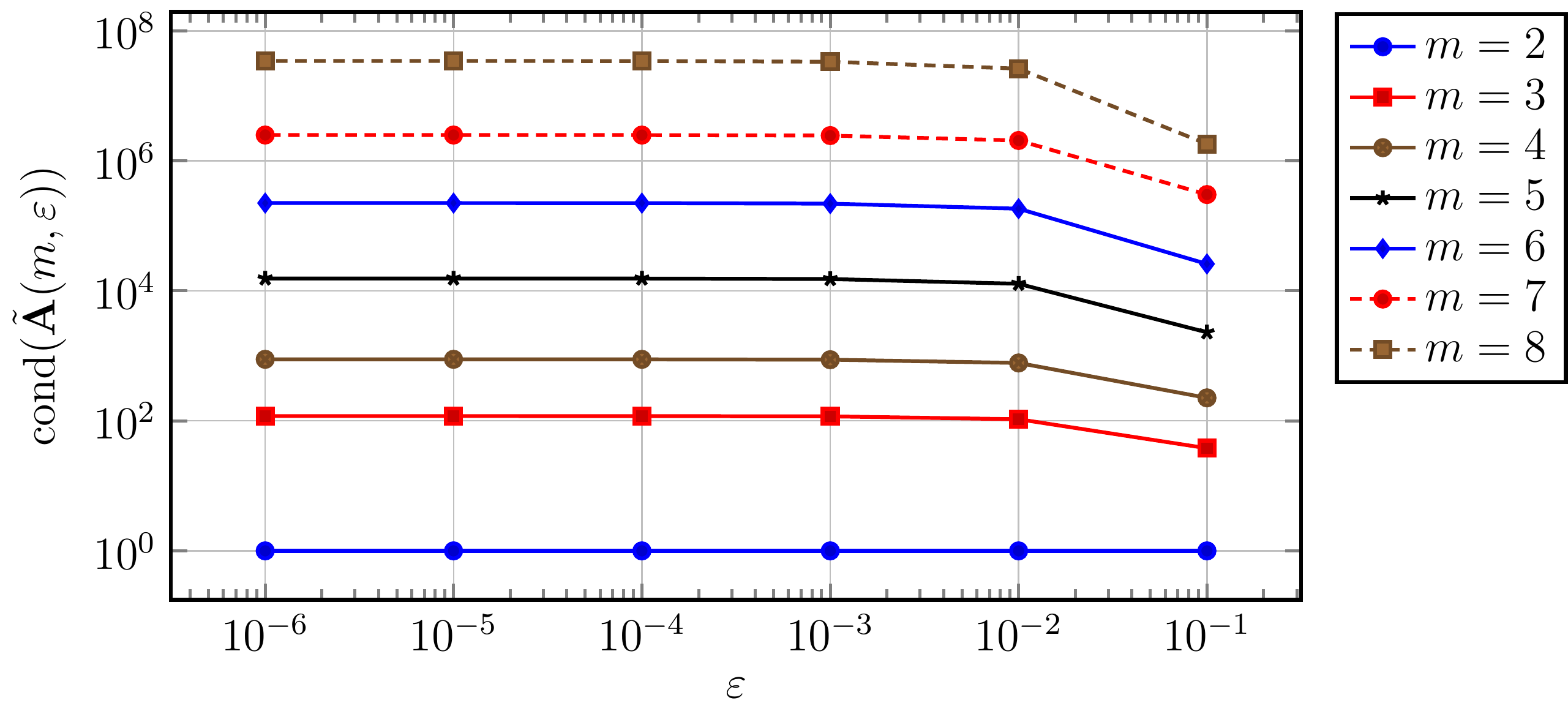}
    \caption{Condition numbers of $\tilde{\mathbf{A}}(m,\varepsilon)$ for different degrees $m$ and different values for $\varepsilon$.}
    \label{fig:cond_numbers_A}
\end{figure}

For demonstrating this desirable feature, consider the case where the interface $\Gamma$ is a circle of radius $1$ and the interface element is $K(\varepsilon)=\frac{1}{\sqrt{2}}-\varepsilon +[0,\frac{1}{2}]^2$ for some $\varepsilon>0$.
Let $\mathbf{A}(m,\varepsilon)$ be the matrix described in \eqref{eqn:def_A} for constructing IFE basis functions on the element $K(\varepsilon)$. We use the basis
$\{p_i(\xi)\}_{i=0}^m = \{L_i(\xi)\}_{i=0}^m$ in \eqref{eq:Q^m_basis_general},
which makes the matrix $\mathbf{A}(m,\varepsilon)$ lower triangular as discussed in the previous section.
Under the assumption that forward substitution is used to solve \eqref{eqn:compact_ife_system} because $\mathbf{A}(m,\varepsilon)$ is lower triangular, the condition of \eqref{eqn:compact_ife_system} can be described by the preconditioned matrix
$\tilde{\mathbf{A}}(m,\varepsilon) = \mathbf{J}^{-1}\mathbf{A}(m,\varepsilon)$ in which $\mathbf{J}=\operatorname{diag}(\mathbf{A}(m,\varepsilon))$ is the Jacobi preconditioner. \autoref{fig:cond_numbers_A} presents data for the condition number of $\tilde{\mathbf{A}}(m,\varepsilon)$ for polynomial's degree $m = 2, 3, \dots, 8$ and for a sequence of $\varepsilon$ values. We note that
as $\varepsilon$ becomes smaller, the sub-element of $K(\varepsilon)$ inside the circular interface $\Gamma$ becomes smaller, and it is extremely small when $\varepsilon = 10^{-6}$. Nevertheless, we observe that the condition number of $\tilde{\mathbf{A}}(m,\varepsilon)$ does not grow as $\varepsilon\to 0$. Similar behavior is observed for the condition number of
$\mathbf{A}(m,\varepsilon)$ when we use the standard basis $\{p_i(\xi)\}_{i=0}^m = \{\xi^i\}_{i = 0}^m$. In this case, since $\mathbf{A}(m,\varepsilon)$ is block lower triangular, the preconditioner $\mathbf{J}$ should be the block diagonal matrix
formed according to $\mathbf{A}(m,\varepsilon)$.

\commentout{
        As \autoref{fig:cond_numbers_A} shows, the condition number appears to approach a positive number that depends on $m$ as $\varepsilon\to 0$. This is a consequence of our choice of the domain of integration in \eqref{eqn:extended_laplacian_IFE_cond}, where we integrate over $\hGamma$ instead of integrating over $\hGamma\cap\hK$. As $\varepsilon\to 0$, we have $|\hGamma\cap\hK|\to 0$. On other hand, $|\hGamma|\to c_0>0$.
}

{
For a summarization, we present a pseudocode in \autoref{algo:construction_frenet_ife} below for constructing a basis
of the local IFE space on an interface element $K$.

\begin{center}
    \begin{algorithm}[H]
        \SetAlgoLined
    \KwData{The interface element $K$, the coefficients $\beta^{\pm}$ and a  parametrization $\g$ of the interface $\Gamma$}
    \KwResult{A basis $\{\phi_i\}_{i=1}^{(m+1)^2}$ of the local IFE space on element $K$}

    Choose $\{\hphi_i\}_{i=m+2}^{(m+1)^2}$ as an enumeration of  $\{\xi^l\eta^j\}_{l=0,j=1}^m$, see \autoref{lem:partial_basis}

    Initialize matrix $\mathbf{A}$ of size  $(m^2-1)\times (m^2-1)$ and matrix $\mathbf{B}$ of size  $(m^2-1)\times (m+1)$\\

    \For{$j=0:m-2$}{
        $\mathbf{A}((m+1)j+1:(m+1)(j+1),:)=A^{(j)}$, \qquad see \eqref{eqn:def_A}

        $\mathbf{B}((m+1)j+1:(m+1)(j+1),:)=[\b^{(j)}(0),\b^{(j)}(1),\dots ,\b^{(j)}(m)]$, \qquad see \eqref{eqn:def_A}

    }

     Solve the system $\mathbf{A}\mathbf{C}= \frac{\beta^--\beta^+}{\beta^+}\mathbf{B}$ \\

    \For{$i=1:m+1$}{
        Let $\c^{(i-1)}=\mathbf{C}(:,i)$

        Set $\hphi_i(\eta,\xi)=p_{i-1}(\xi)$ if $\eta<0$ and by \eqref{eqn:phip_j_2_single_sum} if $\eta>0$

    }

    \Return{$\phi_i(\x) = \hphi_i\circ R_\Gamma(\x), 1 \leq i \leq (m+1)^2$}, \quad see \eqref{eq:IFE_space_span}
\caption{The construction of a basis for the local IFE space on an interface element $K$}
\label{algo:construction_frenet_ife}
    \end{algorithm}
\end{center}













%

}
\section{The approximation capabilities of the IFE space}\label{sec:Approximation}

In this section, we investigate the approximation capabilities of the proposed IFE space on an interface element $K\in \T^i_h$. Since our IFE shape functions are constructed initially on $\hK_F$, it makes sense to obtain error estimates on $\hK_F$ first, then we use the Frenet transformation to obtain similar estimates on $K$. In summary, we start by recalling the classical polynomial approximation results on $\hK_F^\pm$  for a function $\hu\in \H^{m+1}(\hK_F,\hGamma_{K_F};\beta)$, then we will prove that this piecewise polynomial approximation satisfies the discrete interface conditions \eqref{eqn:extended_laplacian_IFE_cond} approximately. After that, we show that approximating $\hu$ on one subelement $\hK_F^{\pm}$ and then extending the resulting  polynomial to an IFE function in $\hV^m_{\beta}(\hK_F)$ will have an optimal $L^2$ error bound. In particular, if the polynomial approximation is defined on the side with the largest coefficient among $\beta^{+}$ and $\beta^-$, then the $L^2$ error bound will not depend on the coefficients $\beta^\pm$. This is a desirable property since the diffusion coefficients can be extremely large or extremely small depending on the application. Furthermore, the error does not depend on the relative position of the interface to the interface element $K$, this property insures that the existence of the so-called small-cut elements will not affect the approximation capabilities of the IFE space. Our approach here is inspired by \cite{guo_higher_2019}, at least in spirit.

For simplicity, we will use the notation $a\lesssim b$ to denote  $a\le C b$ where $C$ is independent of the mesh size, the relative position of the interface and the diffusion coefficients $\beta^\pm$, and will use $a\simeq b$ to denote the equivalence relation $a\lesssim b$ and $b\lesssim a$. Also, for convenience, we will use $s$ to denote a sign $s\in\{+,-\}$ and $s'$ to denote the dual of $s$, that is, if $s=\pm$, then $s'=\mp$.

In the remainder of this section, we will assume that $\g$ is a regular $C^{m+2}([\xi_s,\xi_e], \Gamma)$ parametrization of $\Gamma$. This implies that $\norm{\g'}\simeq 1$, $\norm{\g^{(k)}}\lesssim 1$ for $0\le k\le m+2$ and  $\kappa \in C^m([\xi_s,\xi_e],\mathbb{R})$. In addition, we will assume that $h \kappa_\Gamma \le \frac{1}{2}$, where $\kappa_{\Gamma}=\max_{\xi \in [\xi_s, \xi_e]} |\kappa(\xi)|$, this insures that $\psi \simeq |\hat{D}P_{\Gamma}|\simeq |DR_\Gamma|\simeq 1$.

\begin{lemma}\label{lem:size_hGamma}
    Let $K$ be an interface element, then
    $$|\hGamma_{K_F}| \simeq h,$$
    where $h=\diam(K)$.
\end{lemma}
\begin{proof}
   First, we prove that $h\lesssim |\hGamma_{K_F}|$. We have
$$h^2\simeq \int_{K} d\x=\int_{\hK} |\hat{D}P_{\Gamma}(\eta,\xi)|
d\eta d\xi \lesssim \int_{\hK}
d\eta d\xi\le  \int_{\hK_F}
d\eta d\xi =|\hGamma_{K_F}|h,
$$
Then, we divide both sides by $h$ to obtain $h\lesssim |\hGamma_{K_F}|$. Next, let $\x_1,\x_2\in K$, then by the multivariate mean value theorem, we have
$$\left|\xi(\x_1) -\xi(\x_2)\right|\le \left(\max_{\x\in K} \norm{\nabla \xi(\x)}\right) \norm{\x_1-\x_2}\le h\left(\max_{\x\in K} \norm{\nabla \xi(\x)}\right),$$
We recall from \eqref{eqn:grad_eta_xi} that $\norm{\nabla \xi} = \norm{\g'(\xi)}^{-1}\psi(\eta,\xi)\lesssim 1$. Therefore,
$$|\hGamma_{K_F}|=\max_{\x_1,\x_2\in K} \left|\xi(\x_1) -\xi(\x_2)\right|\lesssim h.$$

\commentout{
\comment{By construction $\hGamma_{K_F}$ is the image of $K$ under $\xi$. Therefore, $$|\hGamma_{K_F}|= \max_{a,b\in \hGamma_{K_F}}\left|a-b\right|= \max_{\x_1,\x_2\in K} \left|\xi(\x_1) -\xi(\x_2)\right|$$ }
}

\end{proof}

In the lemma above, we have shown that the dimensions of the rectangle $\hK_F$ are comparable to those of $K$. In particular, we have proved that $h\lesssim |\hGamma_{K_F}|$. This corroborates out findings at the end of \autoref{sec:Frenet_space}, where we have observed that the condition number of the local IFE system does not grow as the interface becomes close to an edge or a vertex.

Next, let $s\in\{+,-\}$ and let $\hPi^s:H^{m+1}(\hK_F^s)\to Q^m(\hK_F^s)$ be the $L^2$ projection, that is, for $u^s \in H^{m+1}(\hK_F^s)$, $\hPi^s u^s \in  Q^m(\hK_F^s)$ is the unique solution to
$$\left(\hPi^s u^s,v_h\right)_{\hK_F^s} =\left(u^{s},v_h\right)_{\hK_F^s},\qquad \forall v_h\in Q^m(\hK_F^{s}). $$

It is well known \cite{ciarlet_78} that the $L^2$ projection converges optimally as stated in the following lemma.

\begin{lemma}\label{lem:classical_projection}
    Let $s\in\{+,-\}$ and let $\hu^s\in H^{m+1}(\hK_F^s)$, then
\begin{equation}
    \left|\hPi^s\hu^s-\hu^s\right|_{i,\hK_F^s} \lesssim h^{m+1-i}|\hu^s|_{m+1,\hK_F^s}.
    \label{eqn:L2_projection_error}
\end{equation}
For all integers $0\le i \le m+1$.
\end{lemma}

Now, since the dimensions of $\hK_F^{\pm}$ are comparable to $h$ as shown in \autoref{lem:size_hGamma}, the classical trace inequality holds on $H^{1}(\hK_F^s)$ as expressed in the following lemma.

\begin{lemma}\label{lem:classical_trace}
    Let $s\in\{+,-\}$ and let $\hu^s\in H^1(\hK_F^s)$, then
    \begin{equation}
        \norm{\hu^s}_{\hGamma_{K_F}}\lesssim h^{-1/2}\norm{\hu^s}_{\hK_F^s}+ h^{1/2}|\hu^s|_{1,\hK_F^s}.
        \label{eqn:classical_trace}
    \end{equation}
\end{lemma}

Next, we combine \autoref{lem:classical_projection} and \autoref{lem:classical_trace} to obtain an estimate on the projection error on $\hGamma_{K_F}$.

\begin{corollary}\label{coro:classical_projection}
    Let $s\in\{+,-\}$ and let $\hu^s\in H^{m+1}(\hK_F^s)$, then
\begin{equation}
    \left|\hPi^s\hu^s-\hu^s\right|_{i,\hGamma_{K_F}} \lesssim h^{m+\frac{1}{2}-i}|\hu^s|_{m+1,\hK_F^s},
    \label{eqn:trace_projection_error}
\end{equation}
for all integers $i=0,\dots, m$.
\end{corollary}

\begin{proof}
    Using \eqref{eqn:classical_trace}, we have
    $$\left|\hPi^s\hu^s-\hu^s\right|_{i,\hGamma_{K_F}} \lesssim h^{-1/2}\left|\hPi^s\hu^s-\hu^s\right|_{i,\hK_F^s}+h^{1/2}\left|\hPi^s\hu^s-\hu^s\right|_{i+1,\hK_F^s},  $$
    which is combined with \eqref{eqn:L2_projection_error} to obtain \eqref{eqn:trace_projection_error}.
\end{proof}

Next, we will use \autoref{coro:classical_projection} to obtain some estimates on the Laplacian of the residual $\hPi^s\hu^s-\hu^s$ on $\hGamma_{K_F}$.

\begin{lemma}\label{lem:derivatives_of_L_jump}
      Let $s\in\{+,-\}$ and let $\hu^s\in H^{m+1}(\hK_F^s)$, then
\begin{equation}
    \norm{\frac{\p^j}{\p\eta^j}\L(\hPi^s\hu^s-\hu^s)}_{\hGamma_{K_F}} \lesssim h^{m+\frac{1}{2}-2-j}|\hu^s|_{m+1,\hK_F^s},
    \label{eqn:trace_L_projection_error}
\end{equation}
for all integers $j=0,1,\dots, m-2$.
\end{lemma}

\begin{proof}
    Recall that from our assumptions, we have $\norm{\g'}\simeq 1$ and $\norm{\g'},\norm{\g''},|\kappa|,|\kappa'|\lesssim 1$. Therefore, it follows immediately from \eqref{eq:J_derivatives} that
    \begin{equation}\left|\frac{\p^lJ_i}{\p \eta^l}(0,\xi)\right|
    \lesssim 1,\qquad  \forall \xi\in \hGamma_{K_F},\quad  0\le l\le m-2,\quad
    0\le i\le 2.
    \label{eqn:bounds_on_J_i}
    \end{equation}
    Now, we use \eqref{eq:eta_d_L} and \autoref{coro:classical_projection} to obtain
    $$\norm{\frac{\p^j}{\p\eta^j}\L(\hPi^s\hu^s-\hu^s)}_{\hGamma_{K_F}} \lesssim \norm{\hPi^s\hu^s-\hu^s}_{j+2,\hGamma_{K_F}}\lesssim h^{m+1/2-2-j}|\hu^s|_{m+1,\hK_F^s}.$$
\end{proof}

We now consider a function $\hu$ that satisfies the interface conditions \eqref{eqn:flattened_interface_conditions}, and we estimate the discontinuities across the Frenet fictitious interface $\hGamma_{K_F}$ for its projection $\hPi^\pm \hu^{\pm}$ in the following lemma.

\begin{lemma}
    Let $\hu\in \H^{m+1}(\hK_F,\hGamma_{K_F};\beta)$ be a function that satisfies \eqref{eqn:flattened_interface_conditions} and let $\hPi^{s}\hu^s$ be the $L^2$ projections of $\hu^s$ onto $Q^m(\hK_F^s)$ for $s\in\{+,-\}$. Then

    \begin{subequations}
    \begin{equation}\norm{\hPi^+\hu^+-\hPi^-\hu^-}_{\hGamma_{K_F}}\lesssim h^{m+\frac12}|\hu^+|_{m+1,\hK_F},
    \label{eqn:Pp_Pm_0}
    \end{equation}
    \begin{equation}
    \norm{\beta^+(\hPi^+\hu^+)_\eta-\beta^-(\hPi^-\hu^-)_\eta}_{\hGamma_{K_F}}\lesssim(\beta^++\beta^-) h^{m-\frac12}|\hu^+|_{m+1,\hK_F},
    \label{eqn:Pp_Pm_1}\end{equation}
    \begin{equation}
    \norm{\beta^+\frac{\p^j}{\p\eta^j} \L(\hPi^+\hu^+)-\beta^-\frac{\p^j}{\p\eta^j} \L(\hPi^-\hu^-)}_{\hGamma_{K_F}}\lesssim (\beta^++\beta^-) h^{m-\frac32-j}|\hu^+|_{m+1,\hK_F},
    \label{eqn:Pp_Pm_2}
    \end{equation}

   for $j=0,1,\dots m-2$.
   \end{subequations}
\end{lemma}

\begin{proof}
    The first  estimate follows from applying the triangle inequality to \autoref{coro:classical_projection}:
    \begin{align*}
    \norm{\hPi^+\hu^+-\hPi^-\hu^-}_{\hGamma_{K_F}}&\le
    \norm{\hPi^+\hu^+-\hu^+}_{\hGamma_{K_F}}+\norm{\hPi^-\hu^--\hu^-}_{\hGamma_{K_F}}
    +\underbrace{\norm{\hu^+-\hu^-}_{\hGamma_{K_F}}}_{=0}\\
    &\lesssim h^{m+\frac12} |\hu^+|_{m+1,\hK_F^+} +h^{m+\frac12} |\hu^-|_{m+1,\hK_F^-} ,
    \end{align*}
   where we used $\hu^+=\hu^-$ on $\hGamma_{K_F}$. We can prove the second estimate similarly:
    \begin{align*}
    \norm{\beta^+ (\hPi^+\hu^+)_{\eta}-\beta^- (\hPi^-\hu^-)_{\eta}}_{\hGamma_{K_F}}&\le
    \beta^+ \norm{(\hPi^+\hu^+- \hu^+)_{\eta}}_{\hGamma_{K_F}}+\beta^- \norm{(\hPi^-\hu^+ - \hu^-)_{\eta}}_{\hGamma_{K_F}}\\
&~~~~+\underbrace{\norm{\beta^+\hu^+_{\eta}-\beta^- \hu^-_{\eta}}_{\hGamma_{K_F}}}_{=0}\\
    &\lesssim\beta^+ h^{m-\frac12} |\hu^+|_{m+1,\hK_F^+} +\beta^- h^{m-\frac12} |\hu^-|_{m+1,\hK_F^-}.
    \end{align*}
    The proof of the last estimate is similar and uses \autoref{lem:derivatives_of_L_jump}.
\end{proof}

We have shown in \autoref{lem:existence_uniqueness_extension} that for $s\in \{+,-\}$ and a polynomial $q^s\in Q^m(\hK_F^s)$, there is a unique polynomial $q^{s'} \in Q^m(\hK_F^{s'})$, where $s'$ is the dual of $s$, such that $q=q^+\chi_{\hK_F^+} +q^- \chi_{\hK_F^{-}}\in \hV^m_{\beta}(\hK_F) $. For simplicity, we will use $(q^-,q^+)$ to denote $q\in \hV^m_{\beta}(\hK_F)$. Next, we define $\E^{s} :Q^m(\hK_F^{s'})\to Q^m(\hK_F^{s})$ to be the unique mappings such that
$$(q^-,\E^+(q^-))\in  \hV^m_{\beta}(\hK_F),\quad  (\E^-(q^+),q^+)\in \hV^m_{\beta}(\hK_F),\qquad
\forall q^-\in Q^m(\hK_F^-)\text{ and } \forall q^+\in Q^m(\hK_F^+). $$
According to \autoref{lem:existence_uniqueness_extension}, the mappings $\E^s$ for $s\in\{+,-\}$ are well-defined. Furthermore, they are linear. The next lemma shows that $\hPi^s\hu^s$ and $\E^s(\hPi^{s'}\hu^{s'})$ are close to each other on $\hGamma_{K_F}$.

\begin{lemma}\label{lem:E_extensions_error}
      Let $\hu\in \H^{m+1}(\hK_F,\hGamma_{K_F};\beta)$ be a function that satisfies \eqref{eqn:flattened_interface_conditions}, then
      \begin{equation}\norm{\frac{\p^j}{\p \eta^j}\left(\hPi^s\hu^s-\E^s(\hPi^{s'}\hu^{s'})\right)}_{\hGamma_{K_F}}  \lesssim \left(1+\frac{\beta^{s'}}{\beta^s}\right) h^{m+\frac12 -j}|\hu|_{m+1,\hK_F},
      \label{eqn:E_extension_error}
      \end{equation}
      for all $j=0,2,\dots, m$ and $s\in\{+,-\}$.
\end{lemma}

\begin{proof}
     Without loss of generality, we will only consider the case $s=+$. To prove \eqref{eqn:E_extension_error} for $j=0$, we recall that $\E^+(\hPi^-\hu^-)= \hPi^-\hu^-$ on $\hGamma_{K_F}$, Therefore, by \eqref{eqn:Pp_Pm_0}, we have

    $$\norm{\hPi^+\hu^+-\E^+(\hPi^{-}\hu^{-})}_{\hGamma_{K_F}}
    =\norm{\hPi^+\hu^+-\hPi^{-}\hu^{-}}_{\hGamma_{K_F}}
    \lesssim h^{m+\frac12 } |\hu|_{m+1,\hK_F}.$$
    Similarly, we have $\beta^+\frac{\p}{\p \eta} \E^+(\hPi^-\hu^-)= \beta^-\frac{\p}{\p \eta}(\hPi^-\hu^-)$ on $\hGamma_{K_F}$. Hence, we can use \eqref{eqn:Pp_Pm_1} to obtain
    \begin{align}
    \norm{\frac{\p}{\p \eta}\left(\beta^+ \hPi^+\hu^+-\beta^+\E^+(\hPi^-\hu^{-})\right)}_{\hGamma_{K_F}}
    &=\norm{\frac{\p}{\p \eta}\left(\beta^+ \hPi^+\hu^+-\beta^-\hPi^-\hu^{-}\right)}_{\hGamma_{K_F}} \notag \\
    &\lesssim(\beta^++\beta^-) h^{m+\frac12 -1}|\hu|_{m+1,\hK_F}, \label{eqn:E_extension_j=1}
    \end{align}
    which implies \eqref{eqn:E_extension_error} for $j=1$. Next, we prove \eqref{eqn:E_extension_error} for $2\le j\le m$ using strong induction. We assume that \eqref{eqn:E_extension_error} holds for $j=0,1,\dots,i_0$, where $1\le i_0\le m$, then we prove it for $j=i_0+1$.

    To avoid lengthy expressions, we will use $w$ to denote $\hPi^+\hu^+-\E^+(\hPi^{-}\hu^{-})$. Hence,  our induction hypothesis can be written as
\begin{equation}
    \norm{\frac{\p^j}{\p\eta^j} w}_{\hGamma_{K_F}} \lesssim
     \left(1+\frac{\beta^{-}}{\beta^+}\right) h^{m+\frac12 -j}|\hu|_{m+1,\hK_F},\qquad j=0,1,\dots,i_0.
     \label{eqn:w_induction_hyp}
\end{equation}
Since $w|_{\hGamma_{K_F}}\in \P^m(\hGamma_{K_F})$, we can use the classical 1D inverse inequality
    on polynomials \cite{ciarlet_78} to obtain
    \begin{equation}\norm{\frac{\p^{i+j}}{\p\xi^i\p\eta^j}w}_{\hGamma_{K_F}}\lesssim
    \left(1+\frac{\beta^{-}}{\beta^+}\right) h^{m+\frac12 -j-i}|\hu|_{m+1,\hK_F},\qquad j=0,1,\dots,i_0,\quad  1\le i\le m.
\label{eqn:w_1d_inverse}
  \end{equation}
Now, we recall that \eqref{eqn:extended_laplacian_IFE_cond} with $j=i_0-1$ yields
$$\int_{\hGamma_{K_F}}\beta^+\frac{\p^{i_0-1}}{\p\eta^{i_0-1}} \L(\E^+(\hPi^-\hu^-)) v=\int_{\hGamma_{K_F}} \beta^-\frac{\p^{i_0-1}}{\p\eta^{i_0-1}} \L(\hPi^-\hu^-) v,\qquad \forall v \in \P^m(\hGamma_{K_F}). $$
Therefore, for any $v \in \P^m(\hGamma_{K_F})$, we have

\begin{align*}
    \int_{\hGamma_{K_F}} \beta^+\frac{\p^{i_0-1}}{\p\eta^{i_0-1}}\L(w)v &=
    \int_{\hGamma_{K_F}} \beta^+\frac{\p^{i_0-1}}{\p\eta^{i_0-1}}\L(\hPi^+\hu^+)v -
    \beta^+\frac{\p^{i_0-1}}{\p\eta^{i_0-1}}\L(\E^+(\hPi^-\hu^-))v \\
    &=\int_{\hGamma_{K_F}} \left(\beta^+\frac{\p^{i_0-1}}{\p\eta^{i_0-1}}\L(\hPi^+\hu^+) -
    \beta^-\frac{\p^{i_0-1}}{\p\eta^{i_0-1}}\L(\hPi^-\hu^-)\right)v
\end{align*}
  We apply Cauchy-Schwarz inequality and use \eqref{eqn:Pp_Pm_2} to obtain
\begin{align}\left|  \int_{\hGamma_{K_F}} \beta^+\frac{\p^{i_0-1}}{\p\eta^{i_0-1}}\L(w)v \right|
&\le \norm{v}_{\hGamma_{K_F}} \norm{\beta^+\frac{\p^{i_0-1}}{\p\eta^{i_0-1}}\L(\hPi^+\hu^+) -
    \frac{\p^{i_0-1}}{\p\eta^{i_0-1}}\beta^-\L(\hPi^-\hu^-)}_{\hGamma_{K_F}}\notag\\
    &\lesssim \norm{v}_{\hGamma_{K_F}} (\beta^++\beta^-)h^{m+1/2-i_0-1} |\hu|_{m+1,\hK_F},\label{eqn:L_w_d_CS}
    \end{align}
Next, we combine \eqref{eqn:bounds_on_J_i}, \eqref{eq:eta_d_L} and  the triangle inequality to obtain
\begin{equation}\norm{\frac{\p^{i_0-1}}{\p\eta^{i_0-1}}\L(w) - \frac{\p^{i_0+1}}{\p\eta^{i_0+1}} w}_{\hGamma_{K_F}} \lesssim \sum_{j=0}^{i_0} \sum_{i=0}^{i_0-j+1}\norm{\frac{\p^{j+i}}{\p\xi^i\p\eta^j}w}_{\hGamma_{K_F}}\lesssim
 \left(1+\frac{\beta^{-}}{\beta^+}\right) h^{m+\frac12 -i_0-1}|\hu|_{m+1,\hK_F},
 \label{eqn:L_is_almost_normal}
 \end{equation}
where the last inequality follows from \eqref{eqn:w_1d_inverse}. Then, by \eqref{eqn:L_w_d_CS}, \eqref{eqn:L_is_almost_normal}  and Cauchy Schwarz inequality

\begin{align*}\left|  \int_{\hGamma_{K_F}} \frac{\p^{i_0+1}w}{\p\eta^{i_0+1}} v \right|
&\le \left|  \int_{\hGamma_{K_F}} \left(\frac{\p^{i_0-1}}{\p\eta^{i_0-1}}\L(w) - \frac{\p^{i_0+1}}{\p\eta^{i_0+1}} w\right) v \right|+\left|  \int_{\hGamma_{K_F}}\frac{ \p^{i_0-1}}{\p\eta^{i_0-1}}\L(w)  v \right|
,\\
&\lesssim \norm{v}_{\hGamma_{K_F}} \left(1+\frac{\beta^{-}}{\beta^+}\right) h^{m+1/2-i_0-1} |\hu|_{m+1,\hK_F}.\end{align*}
Finally, we replace $v$ with $\frac{\p^{i_0+1}w}{\p\eta^{i_0+1}}(0,\cdot)$, to get
$$\norm{\frac{\p^{i_0+1}w}{\p\eta^{i_0+1}} }_{\hGamma_{K_F}}
\lesssim  \left(1+\frac{\beta^{-}}{\beta^+}\right) h^{m+1/2-i_0-1} |\hu|_{m+1,\hK_F},
$$
which proves \eqref{eqn:E_extension_error} for $j=i_0+1$.

\end{proof}

\begin{lemma}\label{lem:E_extensions_L2}
        Let $\hu\in \H^{m+1}(\hK_F,\hGamma_{K_F};\beta)$ be a function that satisfies \eqref{eqn:flattened_interface_conditions} and let $s\in\{+,-\}$, then
      \begin{equation}\left|\hPi^s\hu^s-\E^s(\hPi^{s'}\hu^{s'})\right|_{i,\hK_F^s}  \lesssim \left(1+\frac{\beta^{s'}}{\beta^s}\right)  h^{m+1 -i}|\hu|_{m+1,\hK_F},
      \label{eqns:E_extension_error_L2}
      \end{equation}
      for $i=0,1,\dots,m+1$.
\end{lemma}

\begin{proof}
    As before, we will restrict our proof to the case $s=+$ and we let $w=\hPi^+\hu^+-\E^+(\hPi^-\hu^-)$. Since $w\in Q^m(\hK_F^+)$, we have for any $(\eta,\xi)\in \hK_F^+$
    $$w(\eta,\xi)=
    \sum_{j=0}^m \frac{\eta^j}{j!} \frac{\p^j w}{\p\eta^j} (0,\xi).
    $$
    By taking the $(k,l)$ derivative of $w$, we obtain

    $$\frac{\p^{k+l}w}{\p\eta^k\p\xi^l} (\eta,\xi)=
    \sum_{j=k}^m \frac{\eta^{j-k}}{(j-k)!} \frac{\p^{j+l} w}{\p\eta^j\p\xi^l} (0,\xi).
    $$
     Then, we take the norm and use the triangle inequality on the right hand side

     \begin{align*}\norm{\frac{\p^{k+l}w}{\p\eta^k\p\xi^l} }_{\hK_F^+}
    &\lesssim \sum_{j=0}^m \left(\int_{0}^h \eta^{2(j-k)}d\eta\right)^{1/2}\norm{\frac{\p^{j+l} w}{\p\eta^j\p\xi^l}}_{\hGamma_{K_F}}\\
    &\lesssim
    \sum_{j=0}^m h^{j-k+1/2}\norm{\frac{\p^{j+l} w}{\p\eta^j\p\xi^l}}_{\hGamma_{K_F}} \lesssim
     \sum_{j=0}^m h^{j-k-l+1/2}\norm{\frac{\p^{j} w}{\p\eta^j}}_{\hGamma_{K_F}},
    \end{align*}
where the last inequality follows from the one-dimensional inverse inequality for polynomials \cite{ciarlet_78}. Now, we apply \autoref{lem:E_extensions_error} to the right hand side
$$
\norm{\frac{\p^{k+l}w}{\p\eta^k\p\xi^l} }_{\hK_F^+}
    \lesssim
    \sum_{j=k}^m h^{j-k-l+1/2}\cdot h^{m+1/2-j} \left(1+\frac{\beta^{-}}{\beta^+}\right)  |\hu|_{m+1,\hK_F}
    =\left(1+\frac{\beta^{-}}{\beta^+}\right)  h^{m+1-l-k }
    |\hu|_{m+1,\hK_F}.
    $$
Finally, \eqref{eqn:E_extension_error} follows from the sum of the estimate above over all non-negative integers $k+l=i$.

\end{proof}

At this point, we are ready to state our first result on the approximation capabilities of the proposed IFE space. For simplicity, we will use $\hPi\hu$ to denote the $L^2$ projection of a function $\hu\in L^{2}(\hK_F)$ onto $\hV^m_{\beta}(\hK_F)$.

\begin{theorem}\label{thm:approx_capabilities_hu}
      Let $\hu\in \H^{m+1}(\hK_F,\hGamma_{K_F};\beta)$ be a function that satisfies \eqref{eqn:flattened_interface_conditions}, then
    \begin{equation}
    \left|\hPi\hu -\hu\right|_{i,\hK_F} \lesssim h^{m+1-i}|\hu|_{m+1,\hK_F},
    \label{eqn:approx_capabilities_hu}
    \end{equation}
    for $i=0,1,\dots,m+1$.
\end{theorem}

\begin{proof}
Without loss of generality, we assume that $\beta^+\ge \beta^-$, then by \eqref{eqn:L2_projection_error} and \eqref{eqns:E_extension_error_L2}, we have

$$\norm{\hu^+-\E^+(\hPi^-\hu^{-})}_{\hK_F^+}  \lesssim   h^{m+1 }|\hu|_{m+1,\hK_F}.
$$
Therefore, by choosing $q=(\hPi^-\hu^-,\E^+(\hPi^-\hu^{-}))$ and using \eqref{eqn:L2_projection_error}, we get
 $$
    \norm{q -\hu}_{\hK_F} \lesssim h^{m+1}|\hu|_{m+1,\hK_F}.
    $$
Since $ \norm{\hPi\hu -\hu}_{\hK_F}\le  \norm{q -\hu}_{\hK_F}$, we obtain \eqref{eqn:approx_capabilities_hu}.
\end{proof}

So far, we have investigated the approximation capabilities of the IFE space on $\hK_F$. In the next step, we will use the mappings $P_\Gamma$ and $R_\Gamma$ to derive error bounds
for the $L^2$ projection of a function $u$ that satisfies the jump conditions \eqref{eqn:continuity_cond}, \eqref{eqn:normal_cond} and \eqref{eqn:laplacian_cond}. For that we essentially follow the ideas in \cite{Zlamel_curved_1973}, where an analysis of an isoparametric finite element method was discussed. Nevertheless, we will show the details for the sake of completeness. First, we define the projection operator $\Pi:L^2(K_F)\to \V^m_{\beta}(K_F)$ using $\hPi$ in the following way
\begin{equation}
    \Pi u = (\hPi \hu)\circ R_{\Gamma},\qquad \hu=u\circ P_{\Gamma}\label{eqn:def_Pi}.
\end{equation}

\begin{theorem}\label{thm:approximation_u}
    Let $u\in \H^{m+1}(K_F,\Gamma;\beta)$ be a function that satisfies \eqref{eqn:laplacian_cond}, then
    $$\left|\Pi u -u \right|_{i,K_F} \le h^{m+1-i}\norm{u}_{m+1,K_F},\qquad 0\le i\le m+1.$$

\end{theorem}
\begin{proof}
    We recall that $\g\in C^{m+2}([\xi_s,\xi_e])$ and $h\kappa_\Gamma\le \frac{1}{2}$. This implies that $P_\Gamma,R_\Gamma \in C^{m+1}([\xi_s,\xi_e])$ and

    \begin{equation}\max_{ \hK_F} \norm{\frac{\p^{k+l}P_\Gamma}{\p\eta^k\p\xi^l}} \lesssim 1,\qquad
    \max_{K_F} \norm{\frac{\p^{k+l}R_\Gamma}{\p x^k\p y^l}} \lesssim 1,
    \label{eqn:dP_dR_bound}
    \end{equation}
    where $0\le k,l\le m+1$. Now, we use \eqref{eqn:def_Pi} to obtain
    $$\left|\Pi u -u \right|_{i,K_F} =\left|\left(\hPi \hu -\hu\right)\circ R_\Gamma \right|_{i,K_F}.  $$
    The change of variables $P_\Gamma$, the inequality $|\hat{D}P_\Gamma|\lesssim 1$ and \autoref{thm:approx_capabilities_hu} leads to
    \begin{equation}  \left|\Pi u -u \right|_{i,K_F} \lesssim
    \left|\hPi \hu -\hu \right|_{i,\hK_F}\lesssim
    h^{m+1-i}|\hu|_{m+1,\hK_F}.
    \label{eqn:Pu_u_hu}
    \end{equation}
    Now, to relate the right hand side to $u$, we recall that $\hu=u\circ P_{\Gamma}$. Therefore, applying the chain rule $m+1$ times yields
    $$
    \left|\frac{\p^{m+1}\hu}{\p \eta^j\p \xi^{m+1-j}}(\eta,\xi) \right|
    \le \left(\sum_{l,k=0}^{m+1}
    \max_{ \hK_F} \norm{\frac{\p^{k+l}P_\Gamma}{\p\eta^k\p\xi^l}}\right)\left(\sum_{k,l=0}^{k+l\le m+1}\left|\frac{\p^{k+l} u}{\p x^k\p y^{l}}\left(P_\Gamma (\eta,\xi)\right)\right|\right).
    $$
    Finally, we take the $L^2$ norm of both sides and use the change of variables $R_\Gamma$
    $$
    |\hu|_{m+1,\hK_F} \lesssim \norm{u}_{m+1,K_F},
    $$
    which, when substituted into \eqref{eqn:Pu_u_hu}, completes the proof.
\end{proof}

\section{An immersed discontinuous Galerkin scheme}\label{sec:the_ipdg_method}

In the previous sections, we have shown that the proposed Frenet IFE space $\V^m_{\beta}(K)$ is easy to construct, approximates the solution $u$ optimally, and the IFE functions within it satisfy the interface conditions \eqref{eqn:continuity_cond} and \eqref{eqn:normal_cond} exactly. This motivates us to employ this space to solve the elliptic interface problem \eqref{eqn:standard_problem_statement}. First, we denote by $\V^m_{\beta}(\T_h)$ the global discontinuous Frenet IFE space defined as
\begin{equation}
    \V^m_{\beta}(\T_h)=\left\{
        \phi \mid \phi|_{K}\in Q^m(K) \text{ if } K\in \T_h^n\text{; Otherwise }\phi|_{K}\in \V^m_{\beta}(K)\right\}.
    \label{eqn:global_DG_space}
\end{equation}
By \autoref{lem:conformity}, we have
$\V^m_{\beta}(\T_h,\Gamma)\subset \H^{s}(\T_h,\Gamma;\beta)$ for any $s>\frac{3}{2}$. This allows us to use the weak formulation \eqref{eqn:compact_weak_form} with no additional integrals on the interface $\Gamma$, in contrast to other high order IFE methods such as the one in \cite{guo_higher_2019}. More specifically, our immersed discontinuous Galerkin (IDG) scheme is as follows

\begin{equation}
\text{Find } u_h\in \V^m_{\beta}(\T_h)\text{ such that} \quad
a_h(u_h,v)=L_h(v),\quad \forall v\in \V^m_{\beta}(\T_h),
    \label{eqn:compact_discrete_weak_form}
\end{equation}
where $a_h$ is defined in \eqref{eqn:def_a_h}. Similarly to the standard interior penalty method \cite{riviere_discontinuous_2008}, the discrete formulation \eqref{eqn:compact_discrete_weak_form} leads to a symmetric linear system $\mathbf{S}\mathbf{c}=\mathbf{f}$. The assembly of $\mathbf{S}$ on non-interface elements and their edges is well known. Hence, we proceed to discuss the assembly of $\mathbf{S}$ on interface elements and their edges. First, we will describe how to implement the Frenet map $R_\Gamma$. After that, we will discuss the construction of the local matrices briefly through an example.

\subsection{The implementation of the Frenet transformation}\label{subsec:implementation_frenet}
Consider an interface element $K$ such as the one shown in \autoref{fig:intf_fic_elements}. Let the points $D$ and $E$ be the intersection points of $\Gamma$ with the boundary of the element. Before implementing $R_\Gamma$, it is useful to perform the following steps
\begin{enumerate}
    \item Find the coordinates of $D$ and $E$ using Newton's method and/or bisection.
    \item Find $\xi_D$ and $\xi_E$ such that $D=\g(\xi_D)$ and $E=\g(\xi_E)$. This can be preformed using Newton's method or the gradient descent applied to $\min_{\xi} |\g(\xi)-\x|^2$ where $\x\in\{D,E\}$.  In our numerical experiments, the gradient descent with Barzilai-Borwein step \cite{barzilai_borwein_1988} outperforms Newton's method and can be implemented to solve $\g(\xi)=\x$ as follows:
    \begin{enumerate}
        \item Perform a line search on the domain of $\g$ to find a good initial condition $\xi_0$.
        \item Compute $\x_1$ using one Newton iteration
        $$\xi_1= \xi_0 -\frac{|\g(\xi_0)-\x|^2}{2\g'(\xi_0)\cdot \left(\g(\xi_0)-\x\right)}.$$
        \item Next, use the following gradient descent iteration until $|\g(\xi_n)-\x|<\varepsilon$ for a set tolerance $\varepsilon$.
        $$\xi_{n+1}=\xi_n -2 \gamma_n \g'(\xi_n)\cdot (\g(\xi_n)-\x),$$
        where the step size $\gamma_n$ is defined in \cite{barzilai_borwein_1988}.

    \end{enumerate}
    \end{enumerate}

Now, we are ready to compute $R_\Gamma(\x)$ for a given $\x$. We let $\overline{\xi}=\frac{1}{2}(\xi_D+\xi_E)$ and use $(\eta_0,\xi_0)=(0,\overline{\xi})$ as an initial guess for Newton's iteration to solve for $R_\Gamma(\x) = (\eta,\xi)$ from $P_\Gamma(\xi,\eta)=\x$:
$$(\eta_{n+1},\xi_{n+1})= (\eta_{n},\xi_{n}) - \hat{D}P_\Gamma(\eta_n,\xi_n)^{-1} P_\Gamma(\eta_n,\xi_{n}),$$
where $\hat{D}P_\Gamma$ is given by \eqref{eqn:jacobian_of_PK}. Alternatively, one can use its inverse \eqref{eqn:jacobian_of_RK}.

\subsection{The construction of the local matrices}\label{subsec:local_matrices}

We construct the local immersed DG matrices on an interface element $K$ by using a basis $\{\phi_i\}_{i=1}^{(m+1)^2}$ of the local IFE space $\hV^m_{\beta}(K)$ as described in \autoref{subsec:basis_IFE}.
For the sake of conciseness, we will only describe the construction of the local stiffness matrix $\mathbf{S}^K$, where $$\mathbf{S}^K_{i,j}= \left( \beta(\x)\nabla\phi_i(\x), \nabla\phi_j(\x)\right)_{K},\qquad i,j=0,1,\dots,(m+1)^2.$$

For that,  we need to compute integrals of piecewise functions on $K$ using a numerical quadrature scheme of the form
$$\int_{K^{\pm}} h(\x) \ d\x \approx \sum_{i=1}^N w_i h(\x_i),
$$
where $(w_i,\x_i)$ are  the weights and nodes of the quadrature. These weights and nodes can be obtained from the procedure discussed in \cite{Saye_quad_2015} where $\x\mapsto \eta(\x)$ plays the role of a level set function. Alternatively, we can split the interface element $K$ into triangles or a quadrilaterals, each on one side of the interface with at most one curved side. After that, we generate the nodes and weights on each triangle/quadrilateral and combine them to form the quadrature rule on $K$. For a detailed discussion, see \cite{moon_immersed_2016}. In our numerical experiments, we did not observe any noticeable difference between the two approaches.

Next, we map the quadrature points $\{\x_i\}_{i=1}^N$ to $\{(\eta_i,\xi_i)\}_{i=1}^N$ using  the Frenet transformation $R_\Gamma$. For $p,q\in \V^m_{\beta}(K)$ we consider the typical integral
$$I=\int_{K} \beta(\x)\nabla p(\x)\cdot \nabla q(\x)\ d\x,$$
which, by \eqref{eqn:def_IIFE_physical}, can be written as
$$I=\int_K\beta(\x) \nabla \hat{p}(\eta(\x),\xi(\x)) \cdot \nabla \hat{q}(\eta(\x),\xi(\x)) \ d\x. $$
By \eqref{eqn:grad_u}, we get
$$I=\int_K\hbeta(\eta,\xi)\left( \hat{p}_{\eta}(\eta,\xi)\hat{q}_{\eta}(\eta,\xi)
+|\g(\xi)|^{-2}\psi(\eta,\xi)^2\hat{p}_{\xi}(\eta,\xi)\hat{q}_{\xi}(\eta,\xi)\right)d\x,$$
where the dependence on $\x$ is dropped for convenience. Lastly, we use the quadrature rule to obtain:

$$I\approx \sum_{i=1}^N w_i\hbeta(\eta_i,\xi_i)\left(\hat{p}_{\eta}(\eta_i,\xi_i)\hat{q}_{\eta}(\eta_i,\xi_i)+|\g(\xi_i)|^{-2}\psi(\eta_i,\xi_i)^2\hat{p}_{\xi}(\eta_i,\xi_i)\hat{q}_{\xi}(\eta_i,\xi_i)\right).$$

In order to compute $\mathbf{S}^K_{i,j}$, we  set $p=\phi_i$ and $q=\phi_j$ in the previous expressions. We compute the other integrals in \eqref{eqn:compact_weak_form} in a similar manner, except that the integrals now are defined on  edges, which are easier to handle.






\section{Numerical Examples}\label{sec:numerical_examples}

In this section, we demonstrate numerically that the proposed immersed DG method converges optimally. Following the estimates of the penalty parameters for the classical symmetric interior penalty method in \cite{epshteyn_estimation_2007}, we set $ \sigma=m^2\sigma_0{\beta^e}$ where $\sigma_0=4$ and ${\beta^e}$ is the maximum value of $\beta$ on an edge $e$. In each example, we fix $\beta^-=1$ and vary $\beta^+\in \{10,100,1000\}$ and present the relative $L^2$ error for different degrees $m=1,2,3,4$ in Examples 1 and 2, and degrees $m=1,2,\ldots,8$ in Example 3. 


\noindent \textbf{Example 1: } Here the domain $\O=(-1,1)^2$ and the interface $\Gamma=\{\x\in\O\mid |\x|=r_0\}$ splits $\Omega$ into $\Omega^-=\{\x\in\O\mid |\x|<r_0\}$ and $\O^+=\{\x\in\O\mid |\x|>r_0\}$, where $r_0$ is set to $\frac{1}{\sqrt{3}}$. The function $f$ and $g$ in \eqref{eqn:Model_problem_elliptic} are chosen such that the exact solution is given by

\begin{equation}u(\x)= \begin{cases}
    \ds \frac{1}{\beta^{-}} \cos\left(\pi |\x|^2\right),& \x\in \O^{-},\\
   \ds  \frac{1}{\beta^{+}} \cos\left(\pi |\x|^2\right) + \cos\left(\pi r_0^2\right) \left(\frac{1}{\beta^-} -\frac{1}{\beta^+}\right),& \x\in \O^+.
\end{cases}\label{eqn:circular_u}
\end{equation}

\begin{figure}[ht]
    \centering
    \begin{subfigure}{.32\textwidth}
        \includegraphics[scale=.15]{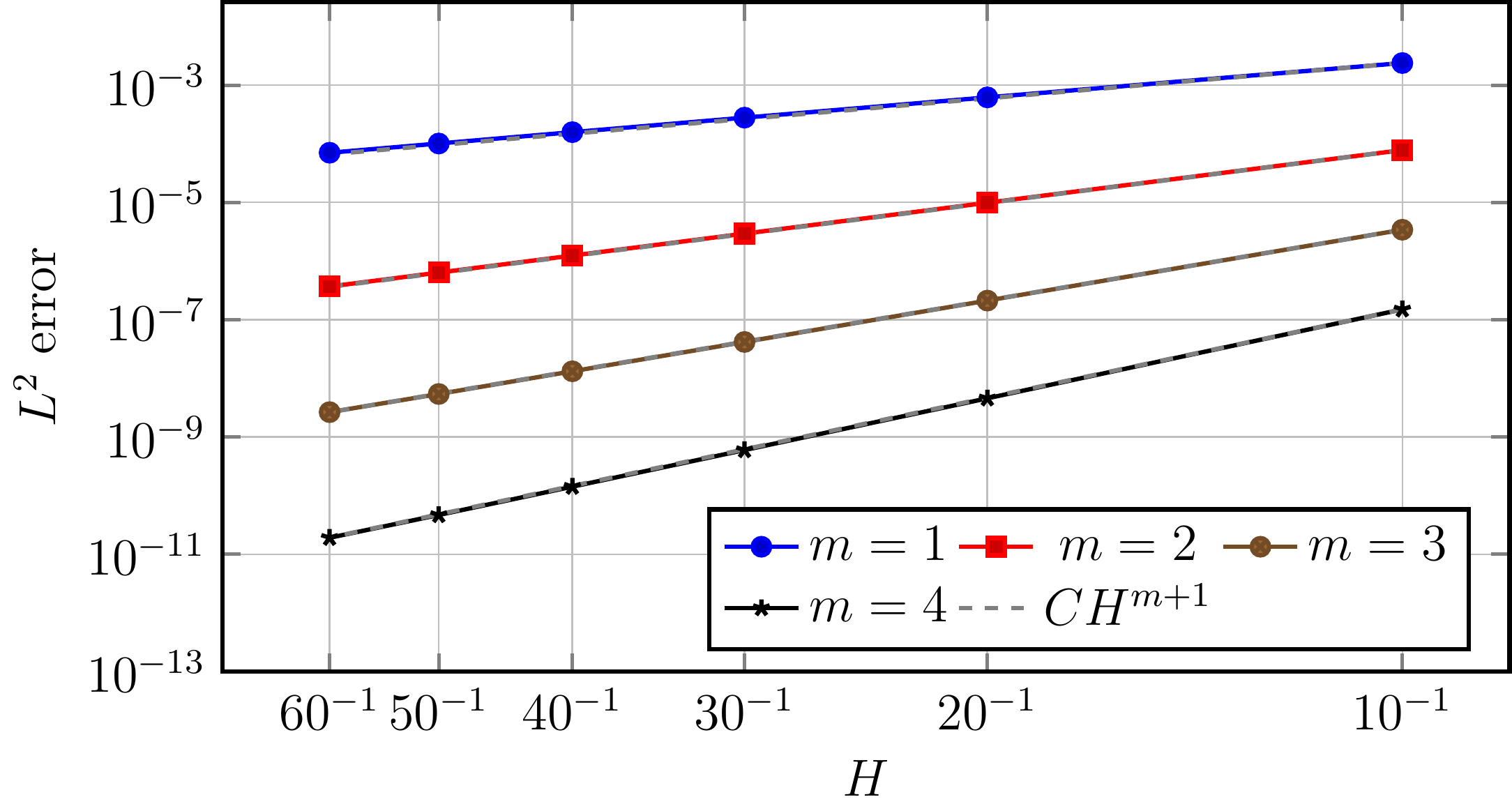}
    \end{subfigure}
     \begin{subfigure}{.32\textwidth}
        \includegraphics[scale=.15]{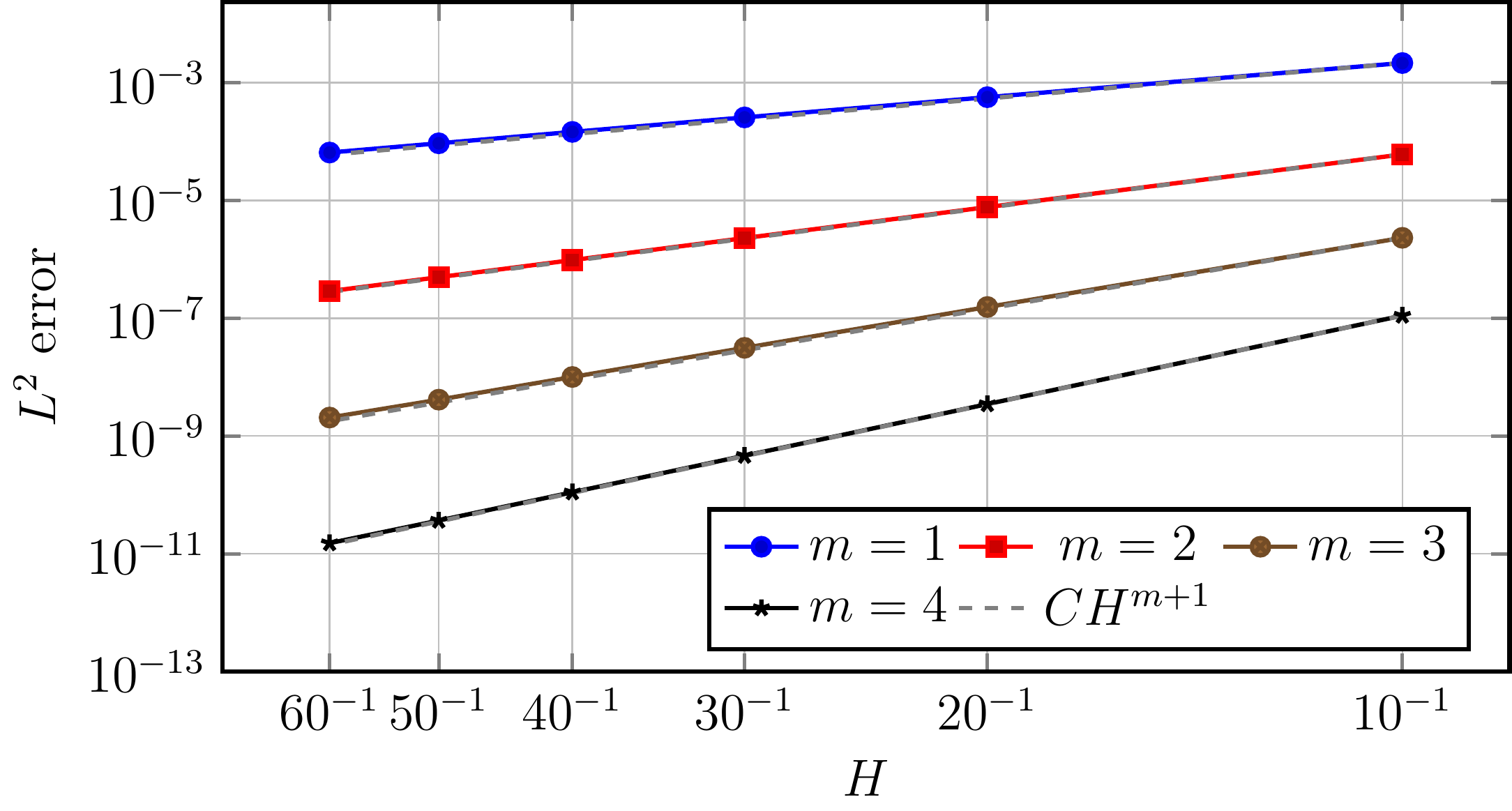}
    \end{subfigure}\begin{subfigure}{.32\textwidth}
        \includegraphics[scale=.15]{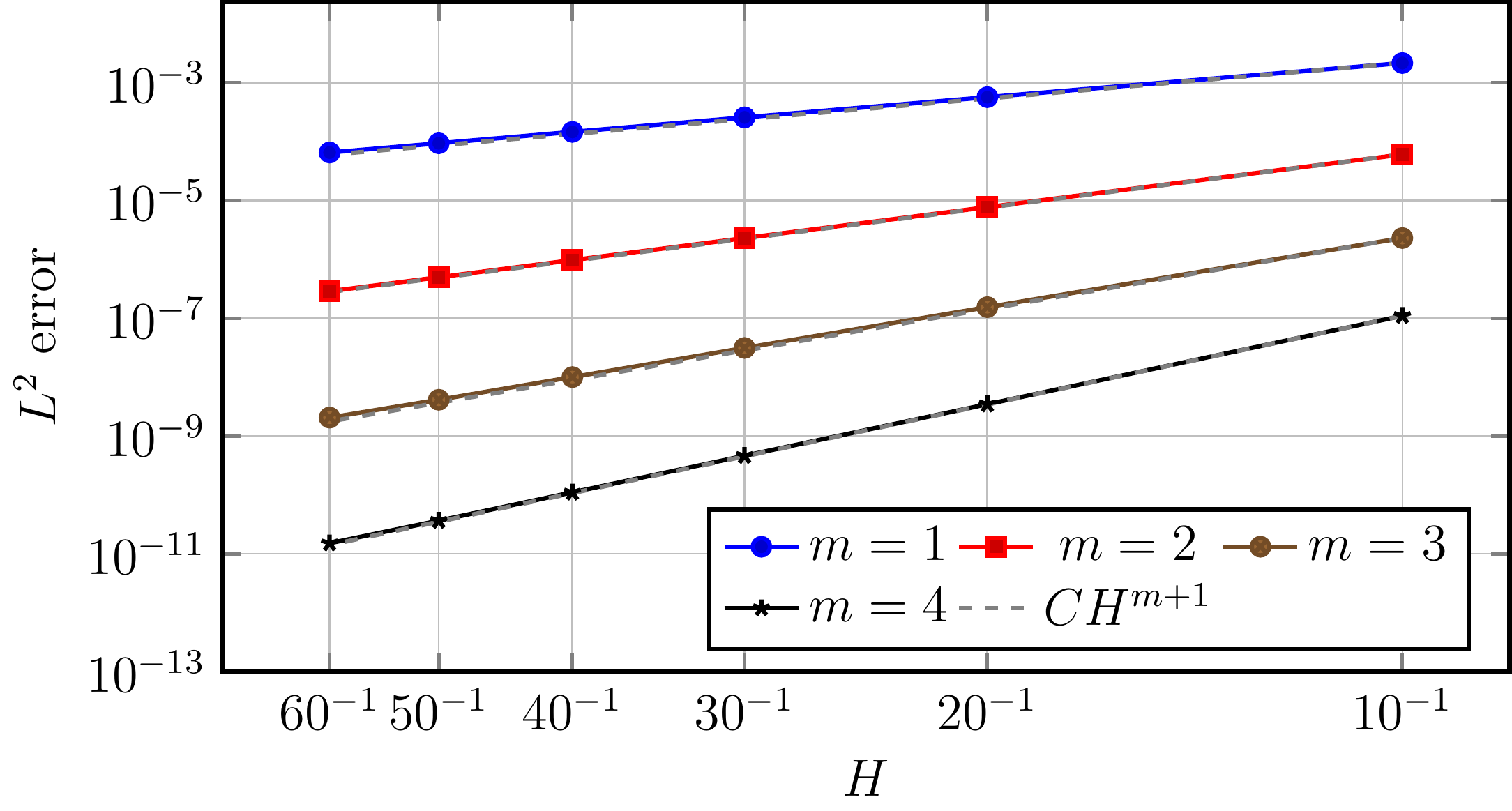}
    \end{subfigure}
    \caption{The relative error of the $L^2$ projection method applied to Example 1 with $\beta^+=10$ ( left), $\beta^+=100$ (middle), $\beta^+=1000$ (right).}
    \label{fig:circular_problem_projection_error}
\end{figure}

First, to illustrate the optimal approximation capability stated in \autoref{thm:approximation_u},  we compute the $L^2$ projection of $u$ onto $\V^m_{\beta}(\T_h)$  on $n\times n$ Cartesian meshes with $n=20,40,\ldots,120$, $m=1,2,3,4$, and plot the resulting relative projection error  versus $H=2/n$ in \autoref{fig:circular_problem_projection_error}. We observe that the $L^2$ projection converges optimally to $u$ under mesh refinement and that it is not affected by the contrast of the coefficients $\beta^{\pm}$.

\begin{figure}[ht]
    \centering
    \begin{subfigure}{.32\textwidth}
        \includegraphics[scale=.15]{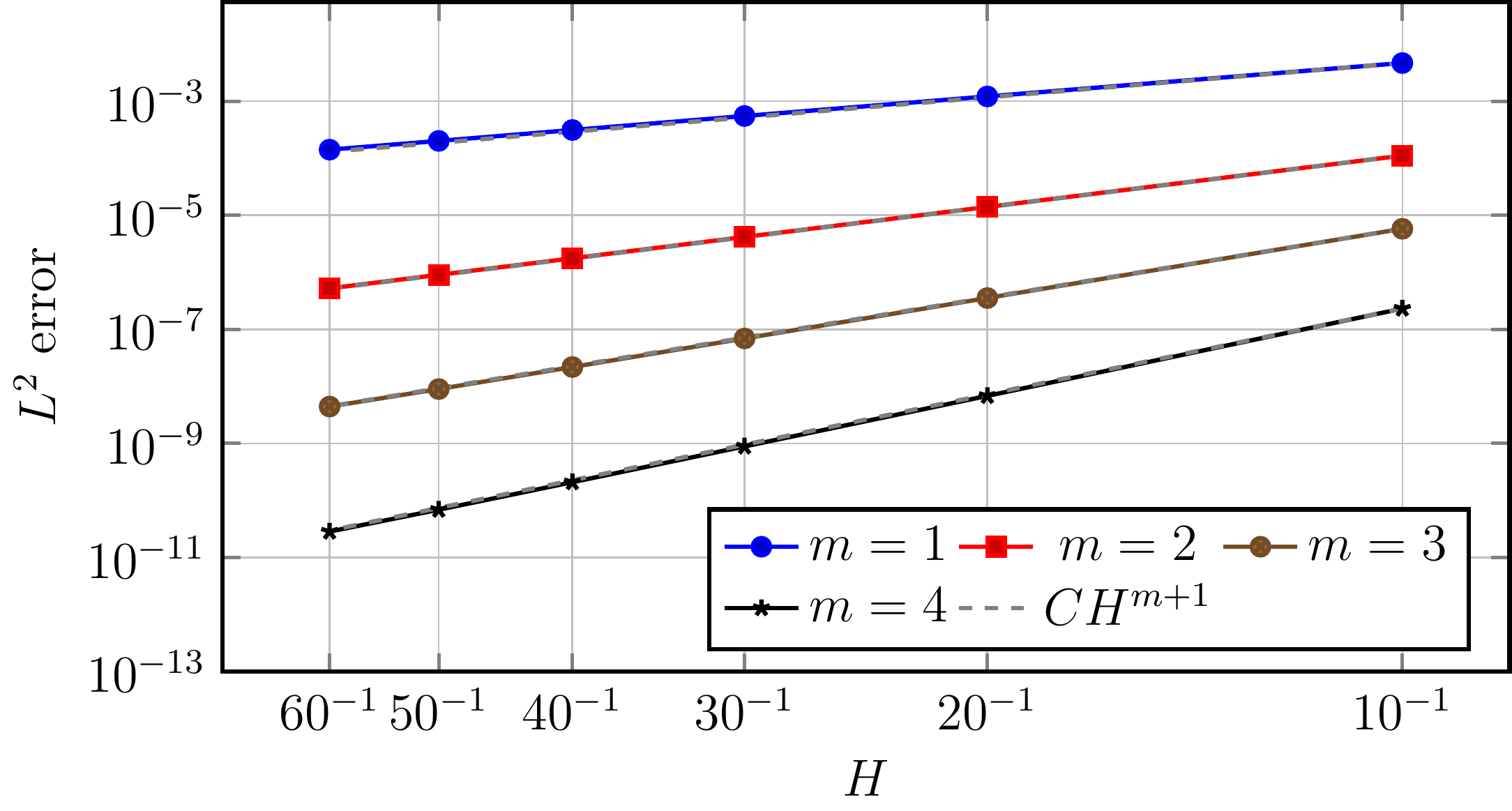}
    \end{subfigure}
     \begin{subfigure}{.32\textwidth}
        \includegraphics[scale=.15]{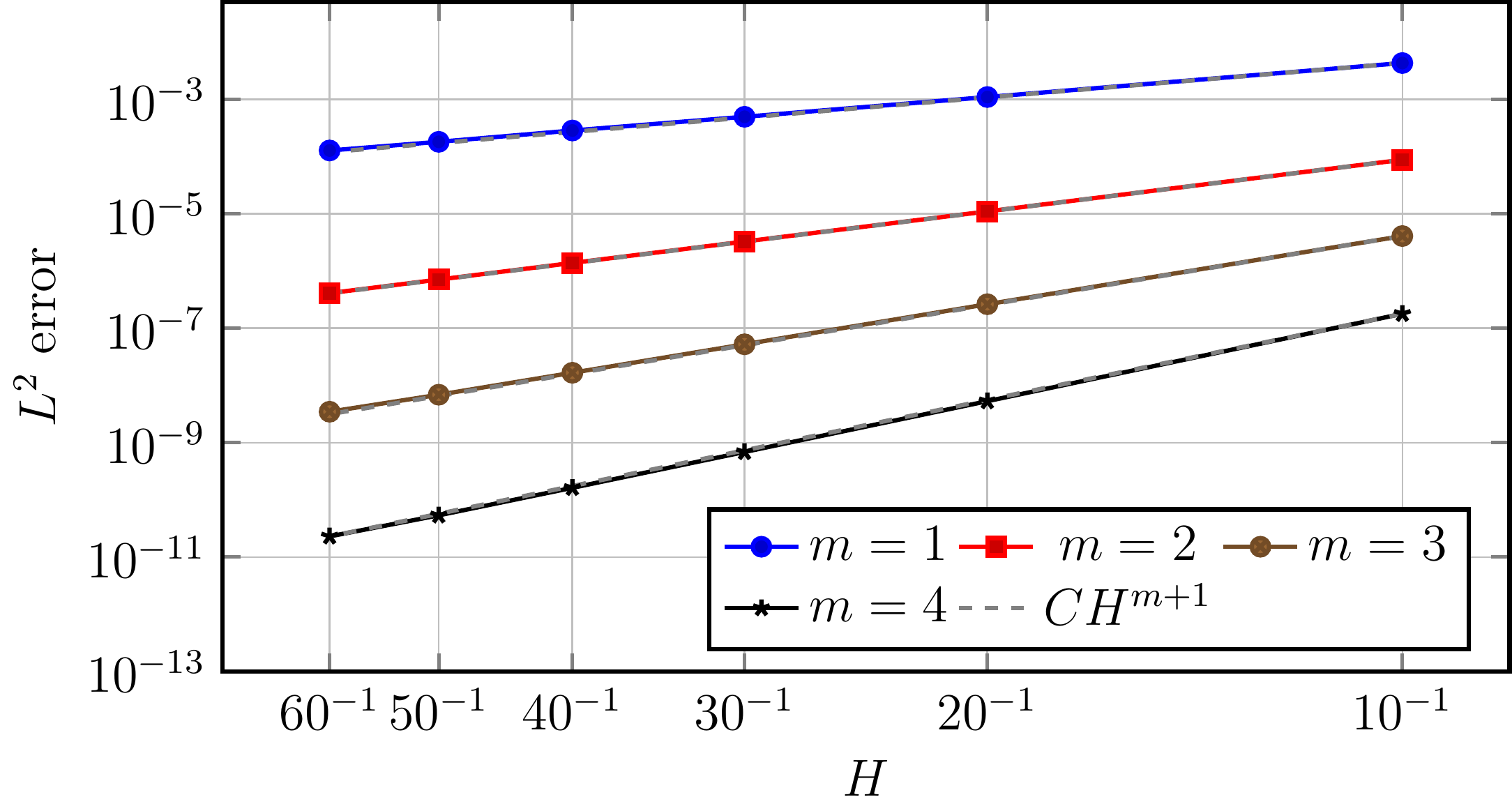}
    \end{subfigure}\begin{subfigure}{.32\textwidth}
        \includegraphics[scale=.15]{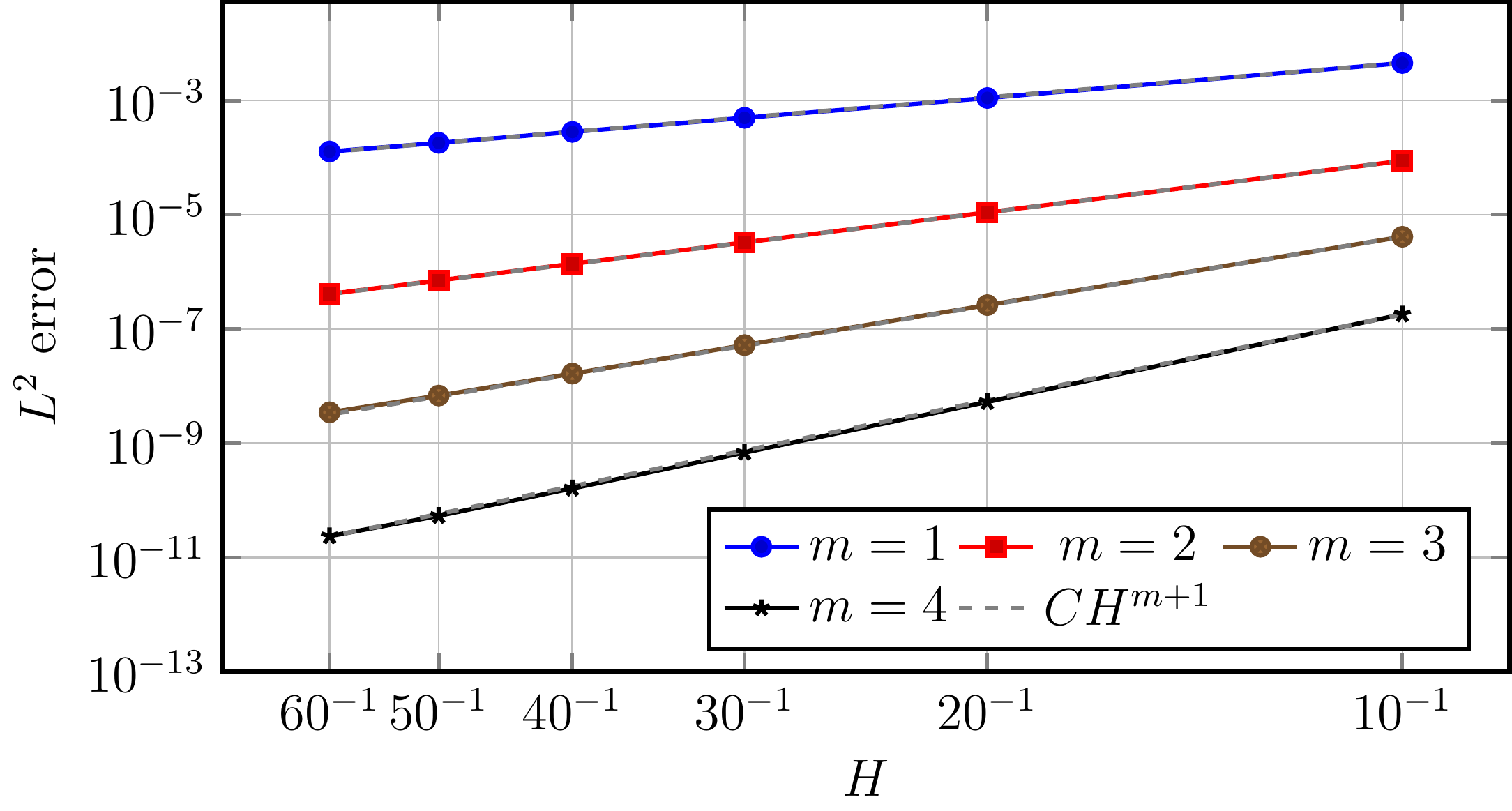}
    \end{subfigure}
    \caption{The relative error of the IDG method applied to Example 1 with $\beta^+=10$ ( left), $\beta^+=100$ (middle), $\beta^+=1000$ (right).}
    \label{fig:circular_problem_method_error}
\end{figure}

We now consider the interface problem \eqref{eqn:standard_problem_statement}  on $\O=(-1,1)^2$ where $f$ and $g$ are selected
such that  the exact solution $u$ is given by \eqref{eqn:circular_u}. We apply our IDG method to solve this problem on Cartesian meshes with $n=20,40,\ldots,120$ and $m=1,2,3,4$ and plot the relative $L^2$ IDG errors  versus $H$
in \autoref{fig:circular_problem_method_error} to observe that the IDG solution converges optimally to $u$ and is not affected by the contrast of the coefficients $\beta^{\pm}$.

\noindent \textbf{Example 2: } In this example, we consider the problem discussed in \cite{adjerid_high_2017}, where $\O=(.6,1.6)\times (.2,1.2)$ and the interface is given by

\begin{equation}\Gamma= \left\{ \x\in \O\mid L(\x):=(x_1^2-x_2^2)^2 -4x_1^2x_2^2+\frac{1}{2}=0\right\}.
\label{eqn:ALR_interface}
\end{equation}
Let $\Omega^{+}= \left\{\x \in \O\mid L(\x)>0\right\}$ and let $\Omega^{-}= \left\{\x \in \O\mid L(\x)<0\right\}$. Next, we define $\tilde{L}$ to be the harmonic conjugate of $L$ given by $\tilde{L}(\x)= 4x_1x_2(x_1^2-x_2^2)$ and we choose $f$ and $g$ such that

$$u(\x) = \frac{1}{\beta(\x)}L(\x) +\tilde{L}(\x) +\frac{1}{\beta(\x)} \tilde{L}(\x)L(\x).$$
Thus, $u$ is a piecewise 7-th degree bi-variate polynomial and the interface $\Gamma$  can be parametrized  by: $$\g(\xi)=\frac{1}{2}\left(\sqrt[4]{2 e^{2 \xi }+1},\sqrt{3 \sqrt{2 e^{2 \xi }+1}-4 e^{\xi }}\right).$$

\begin{figure}[ht]
    \centering
    \begin{subfigure}{.32\textwidth}
        \includegraphics[scale=.15]{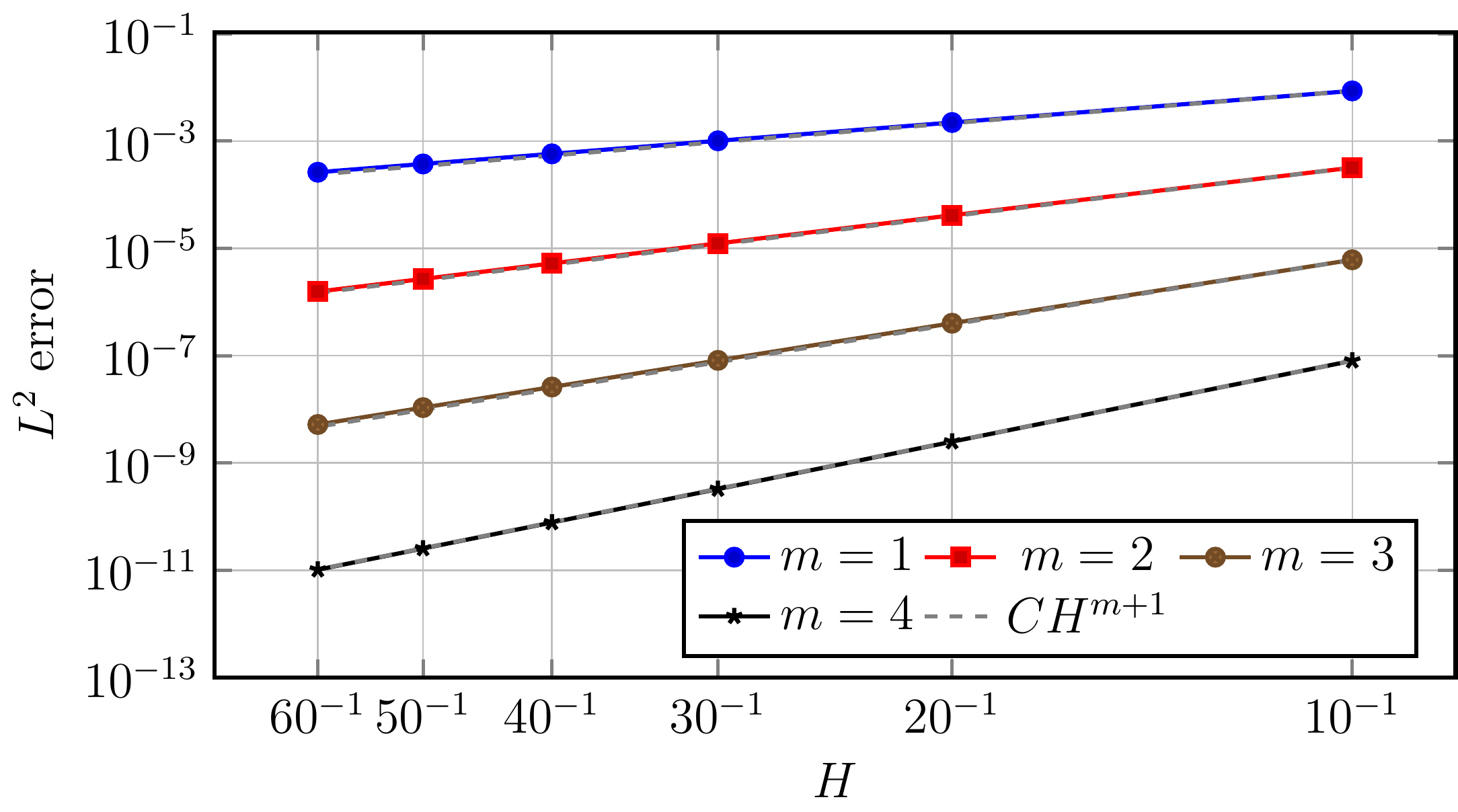}
    \end{subfigure}\begin{subfigure}{.32\textwidth}
        \includegraphics[scale=.15]{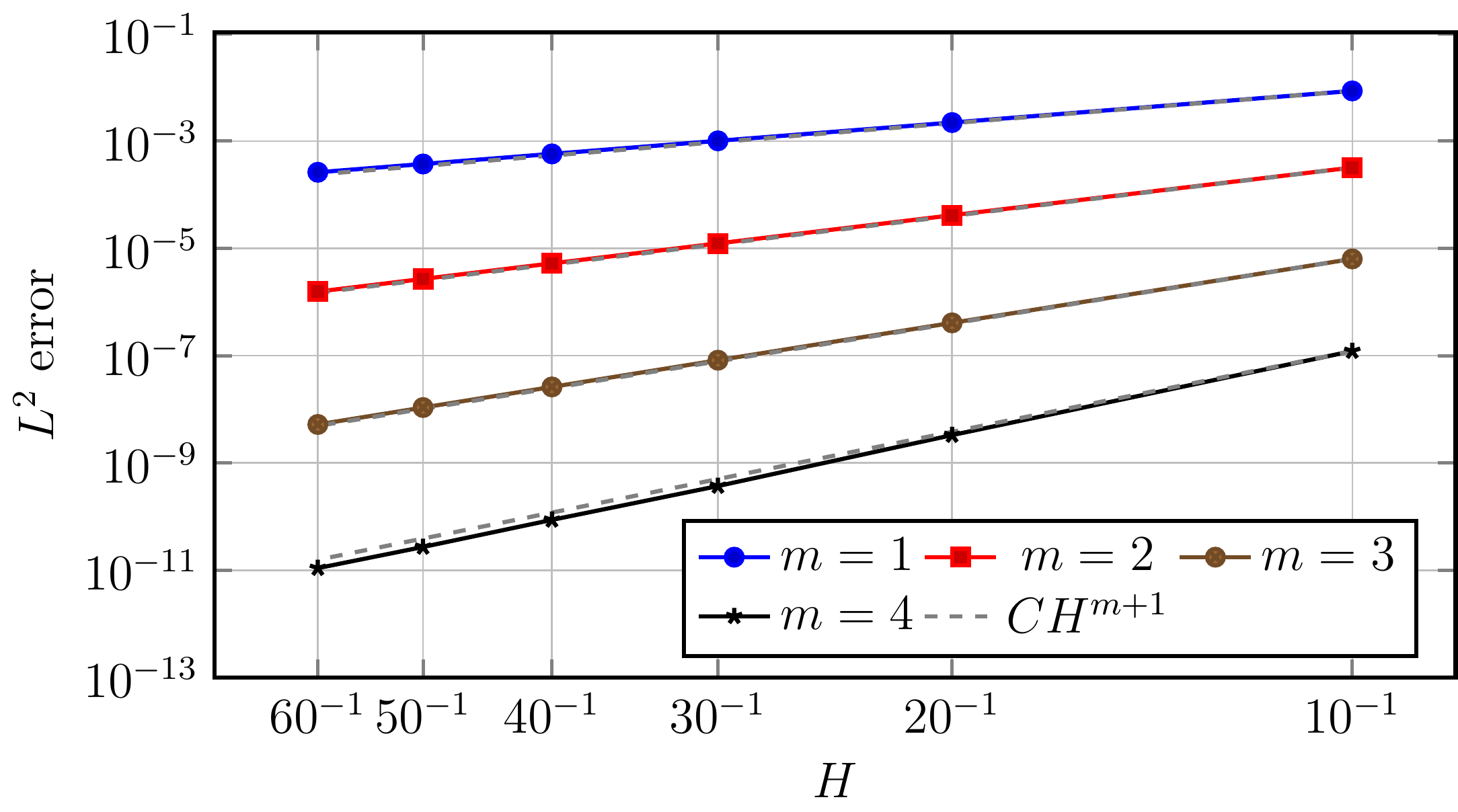}
    \end{subfigure}\begin{subfigure}{.32\textwidth}
        \includegraphics[scale=.15]{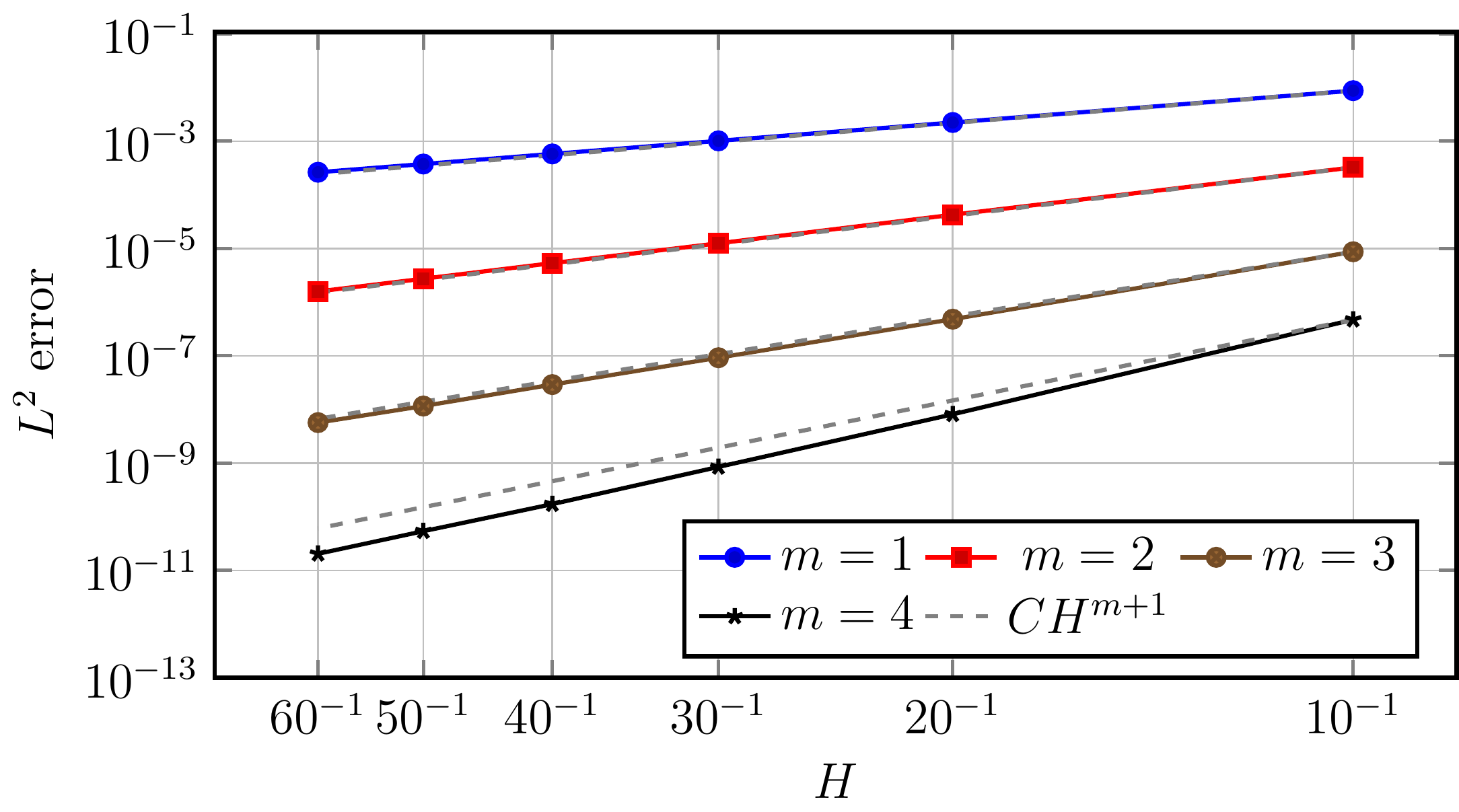}
    \end{subfigure}
    \caption{The relative error of the IDG method applied to Example 2 with $\beta^+=10$ (left), $\beta^+=100$ (center), $\beta^+=1000$ (right).}
    \label{fig:ALR_problem_method_error}
\end{figure}

This example showcases the ability of the proposed IFE space to handle solutions that have non-vanishing tangential derivatives on $\Gamma$ as well as interfaces with non-linear curvature such as \eqref{eqn:ALR_interface}.  We solve this problem on Cartesian meshes with $n=10,20,\ldots,60$, $m=1,2,3,4$, and plot the relative  IDG errors in \autoref{fig:ALR_problem_method_error} versus $H=1/n$ to observe that the method converges at an optimal rate of $H^{m+1}$ under mesh refinement.  \\

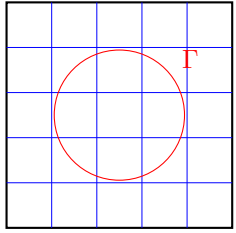
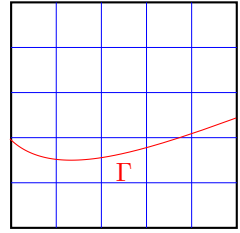
\begin{figure}[ht]
    \centering
    \begin{subfigure}{.4\textwidth}
    \centering
        \begin{tikzpicture}[scale=1.5]
       \draw[thick,black]   (-1,-1)  -- (-1,1) -- (1,1) --(1,-1)   -- cycle;
       \foreach \x in {-.6,-.2,.2,.6 }{
       \draw[very thin,blue] (\x,-1) -- (\x,1);
        \draw[very thin,blue] (-1,\x) -- (1,\x);
       }
       \draw[red] (0,0) circle ({1/sqrt(3)}) ;
       \node[red] at (.63,.5) {$\Gamma$};
    \end{tikzpicture}
    \caption{A $5\times 5$ mesh for Problem $1$}
    \end{subfigure}\begin{subfigure}{.4\textwidth}
    \centering
        \begin{tikzpicture}[scale=3]
       \draw[thick,black]    (.6,.2) -- (.6,1.2) -- (1.6,1.2)--(1.6,.2)   -- cycle;
       \foreach \x in {.2,.4,...,.8}{
       \draw[very thin,blue] (\x+.6,.2) -- (\x+.6,1.2);
        \draw[very thin,blue] (.6,.2+\x) -- (1.6,.2+\x);
       }
       \draw [red,  domain=0:0.379112, samples=60]
 plot ({2^(-.25)*sec(4*\x*180/pi)^(.25)*sin(45+\x*180/pi)}, {2^(-.25)*sec(4*\x*180/pi)^(.25)*sin(45-\x*180/pi)} ) ;
 \node[red] at (1.1,.45){$\Gamma$};
    \end{tikzpicture}
    \caption{A $5\times 5$ mesh for Problem $2$}
    \end{subfigure}

\caption{A $5\times 5$ mesh for solving Problem 1 (left) and Problem 2 (right). }

\label{fig:coarse_meshes}
\end{figure}

\noindent \textbf{Example 3:} This example is for demonstrating the $p$-convergence of the proposed higher degree IDG method, that is, the convergence of the IFE solution with respect to the degree $m$ in the FE/IFE space on a fixed mesh. Here, we refer to the problems in Examples $1$ and  $2$ as Problem 1 and Problem 2, respectively. { For both problems, the IFE solutions are produced on a mesh of size $5\times 5$ on $\Omega$, where $\Omega= [-1,1]^2$ for Problem $1$ and $\Omega= [0.6,1.6]\times [.2,1.2]^2$ for Problem $2$  as shown in \autoref{fig:coarse_meshes}.}

\begin{figure}[ht]
    \centering
    \begin{subfigure}{.5\textwidth}
    \centering
    \includegraphics[scale=.2]{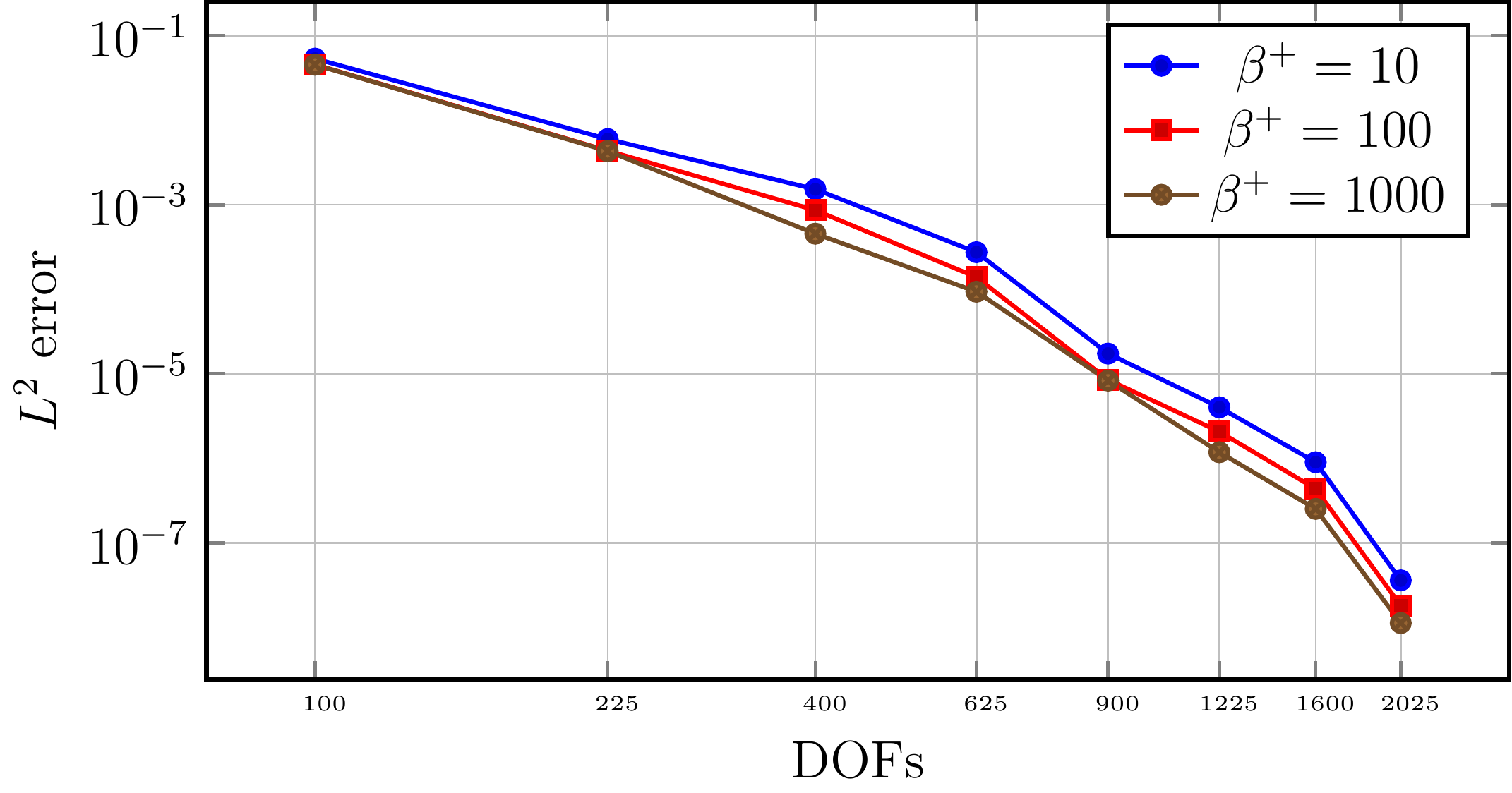}
    \end{subfigure}\begin{subfigure}{.5\textwidth}
        \centering
        \includegraphics[scale=.2]{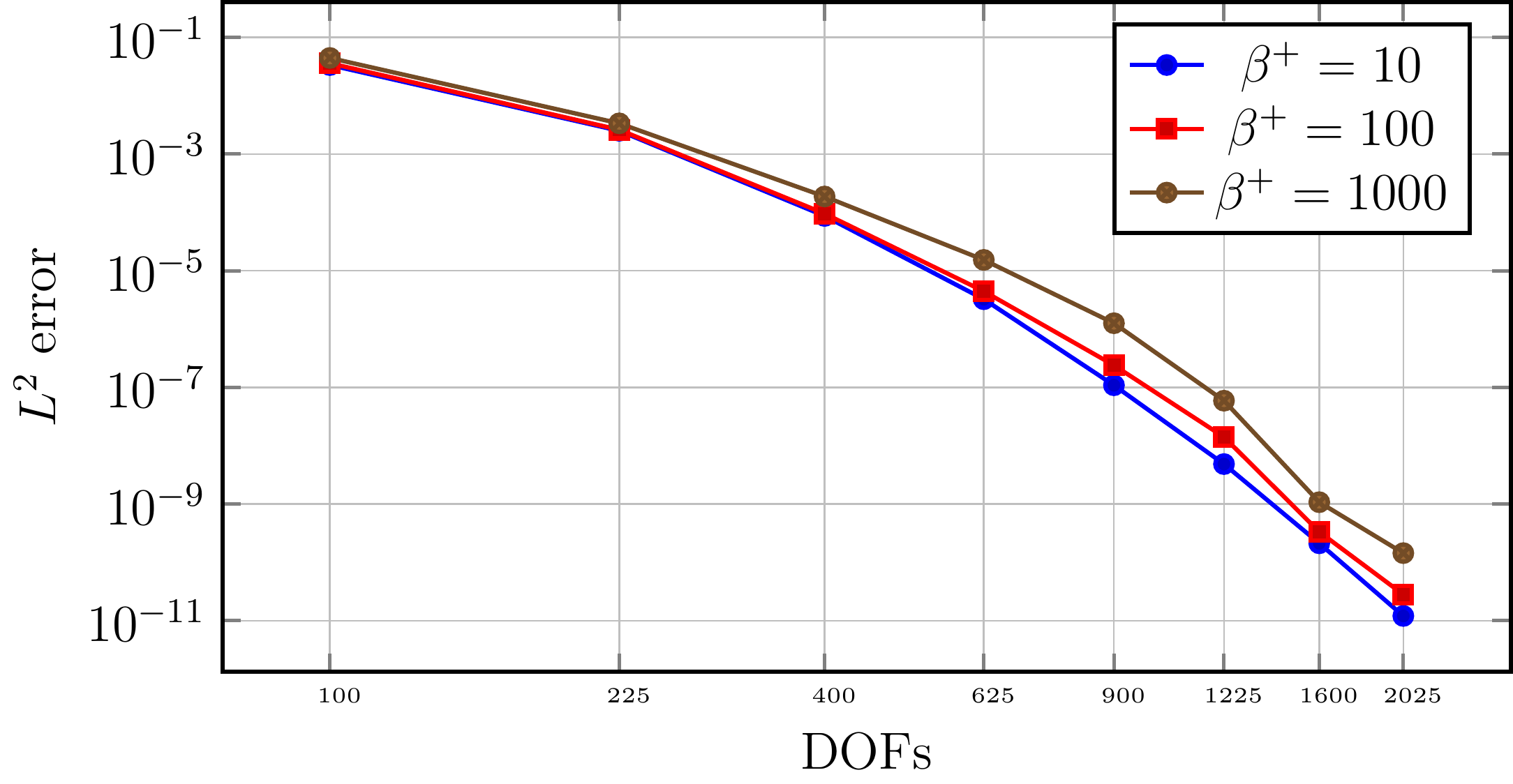}
    \end{subfigure}
    \caption{The relative errors (in log-log plot) of the IDG method applied to Problem 1 (left), and Problem 2 (right) with a fixed mesh and varying degrees $1\le m \le 8$, the $x$-axis represent the total number of degrees of freedom  $25(m+1)^2$.}\label{fig:spectral_error_dofs}.

\end{figure}

 The $L^2$ IDG errors  versus the total number of degrees of freedom shown in \autoref{fig:spectral_error_dofs} in log-log scale demonstrate the exponential convergence of the IDG method with increasing degree $m$ in the presence of a discontinuity of the flux function across the interface. It is worth noting that the interface conditions \eqref{eqn:continuity_cond} and \eqref{eqn:normal_cond} are satisfied exactly by the proposed IFE functions, and, we believe, such a desirable and distinct feature makes it possible for this IDG method to produce very accurate solutions on coarse meshes such as those in \autoref{fig:coarse_meshes}.

In terms of efficiency, this $p$-convergence enables the proposed IDG method to outperform other IFE methods by several orders of magnitude. For instance, when the meshes shown in \autoref{fig:coarse_meshes}
are used with the IFE space $\V^8_{\beta}(\T_h)$, { the total number of degrees of freedom of the proposed method is $2025$}. {Nevertheless, the IFE solutions produced by
the proposed method on such coarse meshes are more accurate than those produced by other IFE methods in the literature, such as
\cite{adjerid_high_2017,guo_group_2019}, requiring more than $100,000$ degrees of freedom.}

%

{Finally, we note that the proposed Frenet framework can be used to construct IFE functions based on bilinear polynomials. In numerical experiments, we observed that the  proposed DG method using the bilinear Frenet IFE functions has an accuracy similar to those bilinear IFE methods in the literature \cite{Adjerid_2023, zhang_bilinear_2015}. In this case, the Frenet framework provides an advantage: The local IFE space $\V^{1}_{\beta}(K)$ can be constructed directly without the need to solve a linear system, which alleviates the ill-conditionning issue associated with small-cut elements. Furthermore, the higher order Frenet space provides a similar error for a much smaller number of degrees of freedom.   }

\section{Conclusion}

A novel high order immersed finite space for the elliptic interface problem is presented. Differential geometry features of the interface are employed in the construction of the IFE functions that satisfy interface jump conditions exactly. Consequently, the proposed IFE space is locally conforming to the Sobolev space for a weak form of the elliptic interface problem. The IFE shape functions on each interface element are easy to construct and the majority of these shape functions can also satisfy the extended jump conditions precisely. The proposed IFE space has been proved to have the optimal approximation capability with respect to the underlying polynomial space. Numerical results demonstrate that the symmetric interior penalty DG method based on this IFE space can solve the elliptic interface problem with the optimal convergence.

In the future, we plan to extend this work in several directions including a rigorous error analysis of the proposed IDG method for elliptic interface problems in two and three dimensions.



\bibliographystyle{siam}
\bibliography{paper_refs}

\begin{thebibliography}{10}

\bibitem{Abate_Tovena_2012}
{\sc M.~Abate and F.~Tovena}, {\em Curves and Surfaces}, Springer-Verlag, Italy, 2012.

\bibitem{gray_2006}
{\sc E.~Abbena, S.~Salamon, and A.~Gray}, {\em Modern {D}ifferential {G}eometry of {C}urves and {S}urfaces with {M}athematica}, Chapman and Hall/CRC, 3~ed., 2006.

\bibitem{Adjerid_2023}
{\sc S.~Adjerid}, {\em A study of high-order immersed finite element spaces by pointwise interface conditions on curved interfaces}, Computers \& Mathematics with Applications, 128 (2022), pp.~331--353.

\bibitem{benromdhane2014}
{\sc S.~Adjerid, M.~Ben-Romdhane, and T.~Lin}, {\em Higher degree immersed finite element methods for second-order elliptic interface problems}, International Journal of Numerical Analysis and Modeling, 11 (2014), pp.~541--566.

\bibitem{adjerid_higher_2018}
\leavevmode\vrule height 2pt depth -1.6pt width 23pt, {\em Higher degree immersed finite element spaces constructed according to the actual interface}, Computers \& Mathematics with Applications, 75 (2018), pp.~1868--1881.

\bibitem{adjerid_high_2017}
{\sc S.~Adjerid, R.~Guo, and T.~Lin}, {\em High degree immersed finite element spaces by a least squares method}, International Journal of Numerical Analysis And Modeling, 14 (2017), pp.~604--625.

\bibitem{barzilai_borwein_1988}
{\sc J.~Barzilai and J.~M. Borwein}, {\em {Two-point step size gradient methods}}, IMA Journal of Numerical Analysis, 8 (1988), pp.~141--148.

\bibitem{benromdhane2015}
{\sc M.~Ben-Romdhane}, {\em Higher-degree immersed finite elements for second-order elliptic interface problems}, PhD thesis, Virginia Tech, 2015.

\bibitem{D.Braess}
{\sc D.~Braess}, {\em Finite {E}lements}, Cambridge University Press, 2~ed., 2001.

\bibitem{burman_cutfem_15}
{\sc E.~Burman, S.~Claus, P.~Hansbo, M.~G. Larson, and A.~Massing}, {\em Cutfem: Discretizing geometry and partial differential equations}, International Journal for Numerical Methods in Engineering, 104 (2015), pp.~472--501.

\bibitem{chorin_numerical_1968}
{\sc A.~J. Chorin}, {\em Numerical solution of the {Navier}-{Stokes} equations}, Mathematics of Computation, 22 (1968), pp.~745--762.

\bibitem{ciarlet_78}
{\sc P.~Ciarlet}, {\em The {F}inite {E}lement {M}ethod for {E}lliptic {P}roblems}, North-Holland, Amsterdam, 1978.

\bibitem{clough_finite_1966}
{\sc R.~Clough and J.~Tocher}, {\em Finite element stiffness matrices for analysis of plates in bending}, in Proceedings of the Conference on Matrix Methods in Structural Mechanics, Wright-Patterson Air Force Base, 1966, pp.~515--545.

\bibitem{epshteyn_estimation_2007}
{\sc Y.~Epshteyn and B.~Rivière}, {\em Estimation of penalty parameters for symmetric interior penalty {Galerkin} methods}, Journal of Computational and Applied Mathematics, 206 (2007), pp.~843--872.

\bibitem{1983Ewingreservoir}
{\sc R.~E. Ewing}, {\em The Mathematics of Reservoir Simulation}, Society for Industrial and Applied Mathematics, 1983.

\bibitem{federer1959curvature}
{\sc H.~Federer}, {\em Curvature measures}, Transactions of the American Mathematical Society, 93 (1959), pp.~418--491.

\bibitem{gersborg-hansen_topology_2006}
{\sc A.~Gersborg-Hansen, M.~P. Bendsøe, and O.~Sigmund}, {\em Topology optimization of heat conduction problems using the finite volume method}, Structural and Multidisciplinary Optimzation, 31 (2006), pp.~251--259.

\bibitem{guo_group_2019}
{\sc R.~Guo and T.~Lin}, {\em A group of immersed finite element spaces for elliptic interface problems}, IMA Journal of Numerical Analsys, 39 (2019), pp.~482--511.

\bibitem{guo_higher_2019}
\leavevmode\vrule height 2pt depth -1.6pt width 23pt, {\em A higher degree immersed finite element method based on a {Cauchy} extension}, SIAM Journal On Numerical Analysis, 57 (2019), pp.~1545--1573.

\bibitem{guo_immersed_2020}
\leavevmode\vrule height 2pt depth -1.6pt width 23pt, {\em An immersed finite element method for elliptic interface problems in three dimensions}, Journal of Computational Physics, 414 (2020), p.~109478.

\bibitem{guzman_2017}
{\sc J.~Guzmán, M.~A. Sánchez, and M.~Sarkis}, {\em A finite element method for high-contrast interface problems with error estimates independent of contrast}, Journal of Scientific Computing, 73 (2017).

\bibitem{he_approximation_2008}
{\sc X.~He, T.~Lin, and Y.~Lin}, {\em Approximation capability of a bilinear immersed finite element space}, Numerical Methods for Partial Differential Equations: An International Journal, 24 (2008), pp.~1265--1300.

\bibitem{guzman2016}
{\sc {J. Guzm\'an}, {M. A. S\'anchez}, and {M. Sarkis}}, {\em Higher-order finite element methods for elliptic problems with interfaces}, ESAIM: Mathematical Modelling and Numerical Analysis, 50 (2016), pp.~1561--1583.

\bibitem{Kuznetsov2006}
{\sc E.~Kuznetsov}, {\em Optimal parametrization in numerical construction of curve}, Journal of The Franklin Institute, 344 (2007), pp.~658--671.

\bibitem{Lane_1940}
{\sc E.~P. Lane}, {\em Metric {D}ifferential {G}eometry of {C}urves and {S}urfaces}, The University of Chicago Press, 1940.

\bibitem{lehrenfeld_high_2016}
{\sc C.~Lehrenfeld}, {\em High order unfitted finite element methods on level set domains using isoparametric mappings}, Computer Methods in Applied Mechanics and Engineering, 300 (2016), pp.~716--733.

\bibitem{leveque_immersed_1994}
{\sc R.~J. LeVeque and Z.~Li}, {\em The immersed interface method for elliptic equations with discontinuous coefficients and singular sources}, SIAM Journal on Numerical Analysis, 31 (1994), pp.~1019--1044.
\newblock Publisher: SIAM.

\bibitem{li_immersed_2004}
{\sc Z.~Li, T.~Lin, Y.~Lin, and C.~Rogers, R.}, {\em An immersed finite element space and its approximation capability}, Numerical Methods for Partial Differential Equations, 20 (2004), pp.~338--367.

\bibitem{li_new_2003}
{\sc Z.~Li, T.~Lin, and X.~Wu}, {\em New {C}artesian grid methods for interface problems using finite element formulation}, Numerische Mathematik, 96 (2003), pp.~61--98.

\bibitem{lin_rectangular_2001}
{\sc T.~Lin, Y.~Lin, R.~C. Rogers, and L.~M. Ryan}, {\em A rectangular immersed finite element method for interface problems}, in Scientific Computing and Applications, vol.~7, Commack, NY, USA, 2001, Nova Science Publishers, Inc., pp.~107--114.

\bibitem{zhang_bilinear_2015}
{\sc T.~Lin, Y.~Lin, and X.~Zhang}, {\em Partially penalized immersed finite element methods for elliptic interface problems}, SIAM Journal on Numerical Analysis, 53 (2015), pp.~1121--1144.

\bibitem{moon_immersed_2016}
{\sc K.~Moon}, {\em Immersed {Discontinuous} {Galerkin} {Methods} for {Acoustic} {Wave} {Propagation} in {Inhomogeneous} {Media}}, PhD thesis, Virginia Tech, May 2016.

\bibitem{ONeill_DiffGeom_2010}
{\sc B.~O'Neill}, {\em Elementary {D}ifferential {G}eometry}, Academic Press, 2~ed., 2006.

\bibitem{Pressley_DiffGeom_2010}
{\sc A.~Pressley}, {\em Elementary {D}ifferential {G}eometry}, Springer Undergraduate Mathemtics Series, Springer, 2~ed., 2010.

\bibitem{riviere_discontinuous_2008}
{\sc B.~Rivière}, {\em Discontinuous {G}alerkin {M}ethods for {S}olving {E}lliptic and {P}arabolic {E}quations}, vol.~FR35 of Frontiers in Applied Mathematics, SIAM, Philadelphia, 2008.

\bibitem{Saye_quad_2015}
{\sc R.~I. Saye}, {\em High-order quadrature methods for implicitly defined surfaces and volumes in hyperrectangles}, SIAM Journal on Scientific Computing, 37 (2015), pp.~A993--A1019.

\bibitem{Kuznetsov2003}
{\sc V.~I. Shalashilin and K.~E. B.}, {\em Parametric Continuation and Optimal Parametrization in Applied Mathematics and Mechanics}, Springer Netherlands, 2023.

\bibitem{Toponogov_2006}
{\sc V.~A. Toponogov}, {\em Differential {G}eometry of {C}urves and {S}urfaces : {A} {C}oncise {G}uide}, Birkh\"{a}user, Boston, 2006.

\bibitem{wang_modeling_2008}
{\sc J.~Wang, X.-M. He, and Y.~Cao}, {\em Modeling electrostatic levitation of dust particles on lunar surface}, IEEE Transactions on Plasma Science, 36 (2008), pp.~2459--2466.

\bibitem{Milliam_Parker_1977}
{\sc R.~S. William and G.~D. Parker}, {\em Elements of {D}ifferential {G}eometry}, Prentice-Hall, 1977.

\bibitem{Zlamel_curved_1973}
{\sc M.~Zl\'{a}mal}, {\em Curved elements in the finite element method. {I}}, SIAM Journal on Numerical Analysis, 10 (1973), pp.~229--240.

\end{thebibliography}

\end{document}